\documentclass{amsart}   	
\usepackage[margin=1in]{geometry}                		
\usepackage{graphicx}				
\usepackage{amssymb,amsfonts,amsmath}
\usepackage{soul}
\usepackage{tikz}
\usepackage{tikz-cd} 

\newtheorem{dummy}{dummy}[section]
\newtheorem{theorem}[dummy]{Theorem}

\newtheorem{conjecture}[dummy]{Conjecture}

\newtheorem{definition}[dummy]{Definition}
\newtheorem{statement}[dummy]{Statement}
\newtheorem{proposal}[dummy]{Proposal}
\theoremstyle{definition}

\newtheorem{example}[dummy]{Example}
\newtheorem{remark}[dummy]{Remark}

\newcommand\cB{\mathcal B}
\newcommand\cC{\mathcal C}
\newcommand\cD{\mathcal D}

\newcommand\cG{\mathcal{G}}

\newcommand\cV{\mathcal V}

\newcommand\cX{\mathcal X}

\newcommand\bC{\mathbb C}

\newcommand\bN{\mathbb N}

\newcommand\bR{\mathbb R}
\newcommand\bZ{\mathbb Z}

\newcommand\rB{\mathrm B}
\newcommand\rT{\mathrm T}
\newcommand\rO{\mathrm O}
\newcommand\PL{\mathrm{PL}}
\newcommand\SO{\mathrm{SO}}

\newcommand\rN{\mathrm N}

\usepackage[utf8]{inputenc}
\DeclareFontFamily{U}{min}{}
\DeclareFontShape{U}{min}{m}{n}{<-> udmj30}{}

\newcommand\define[1]{\emph{#1}}
\newcommand\cat[1]{\mathbf{#1}}

\newcommand\op{\mathrm{op}}
\newcommand\id{\mathrm{id}}
\newcommand\fd{\mathrm{fd}}
\newcommand\pt{\mathrm{pt}}
\newcommand\fr{\mathrm{fr}}

\DeclareMathOperator\Aut{Aut}
\newcommand\Cat{\cat{Cat}}
\newcommand\Gpd{\cat{Gpd}}
\newcommand\dagCat{{\dagger}\cat{Cat}}
\newcommand\flagdagCat{\cat{Fl{\dagger}Cat}}
\newcommand\coflagdagCat{\cat{coFl{\dagger}Cat}}

\newcommand\SymRigidCat{\cat{SymRigidCat}}

\newcommand\AdjCat{\cat{AdjCat}}
\newcommand\Bord{\cat{Bord}}
\newcommand\sVect{\cat{sVect}}
\newcommand\sHilb{\cat{sHilb}}
\newcommand\Diff{\mathrm{Diff}}

\usepackage{xcolor}
\usepackage[plainpages=false,hypertexnames=false,pdfpagelabels]{hyperref}
\definecolor{medium-blue}{rgb}{0,0,.8}
\hypersetup{colorlinks, linkcolor={purple}, citecolor={medium-blue}, urlcolor={medium-blue}}
\newcommand{\arxiv}[1]{\href{http://arxiv.org/abs/#1}{\tt arXiv:\nolinkurl{#1}}}
\newcommand{\arXiv}[1]{\href{http://arxiv.org/abs/#1}{\tt arXiv:\nolinkurl{#1}}}

\newcommand{\doi}[1]{\href{http://dx.doi.org/#1}{{\tt DOI:#1}}}

\newcommand\wipedge{volutive}
\newcommand\wipedgenoun{volution}

\newcommand\RENAMEME{$\PL$-dagger tower}

\title{Dagger $n$-categories}
\author[GF, BH, TJF, CK, LM, N, DP, DR, CS, LS, CV]{Giovanni Ferrer, Brett Hungar, Theo Johnson-Freyd, Cameron Krulewski, Lukas M\"{u}ller, Nivedita, David Penneys, David Reutter, Claudia Scheimbauer, Luuk Stehouwer, Chetan Vuppulury}

\begin{document}
\maketitle

\begin{abstract}
Category theory provides a unified language for organizing composable operations in many disciplines. In disciplines where unitarity is fundamental---such as functional analysis, quantum field theory, and quantum logic---this language must also capture adjoints, leading to the notion of dagger categories. Higher category theory, which extends this framework to encode operations between operations, has recently become indispensable in both theoretical physics and pure mathematics.
Finding a higher categorical analogue of a dagger category is therefore key to the foundations of 
quantum field theory.

{In this work,
    we present a coherent definition of \define{dagger $(\infty,n)$-category} in terms of equivariance data trivialized on parts of the category. Our main example is the bordism $(\infty,n)$-category $\Bord_{n}^X$. This allows us to define (fully-local) \define{reflection-positive topological quantum field theories} to be higher dagger functors out of $\Bord_{n}^X$. 
}

\end{abstract}


\section{Introduction}
Category theory is the main language for organizing mathematical objects and the operations between them. 
It has been applied 
widely throughout computer science---and more recently, biology and engineering---to study how operations, resources, and other quantities compose.
In theoretical physics, categories model systems of elementary particles that merge and split.
In many contexts, notably those coming from functional analysis and physics, each operation $f\colon X\to Y$ has an \emph{adjoint} operation $f^\dagger\colon Y \to X$, 
generalizing
the conjugate transpose of a matrix.
Tracking the assignment $f \mapsto f^\dagger$ is vital for analyzing quantum systems; for example, the measurement corresponding to a quantum operation $f \colon X \to X$ takes real values if and only if  $f^\dagger = f$.
\emph{Dagger categories}, which are equipped with such an anti-involution, organize and axiomatize collections of operations with both compositions and adjoints.

Often, the operations between objects in a category are themselves acted upon by higher operations. Examples include procedures that convert one
computer program into another and interfaces between defects in quantum materials.
Axiomatizing these higher-order processes between processes naturally leads to the concept of \emph{higher categories}. 

In many mathematical and physical applications, 
 higher processes also admit adjoints.
This is particularly apparent in quantum field theory and quantum condensed matter. Quantum field theory is based on two foundational precepts: locality and unitarity. The modern approach to encoding locality is through higher categories, whereas the modern approach to encoding unitarity is through dagger categories.
However, until now there has not been an axiomatization that accommodates both adjoints and higher operations. Our goal is to answer the question:
\textbf{What is a dagger higher category?}

To explain why this  is not straightforward to answer, let us recall the traditional definition of dagger category. 
Categories with an (anti-)involution have been studied since the early days of category theory~\cite{MR0126477,MR0141698}. 
The notions of $\rm C^*$- and $\rm W^*$-category in the complex linear setting were introduced in \cite{MR0808930} 
and have been studied extensively in the $\rm C^*$ tensor setting \cite{MR1239440,MR1010160,MR1444286,MR1470857,MR1966524,MR2091457}.
The complete nonlinear definition is due to Selinger \cite{selinger2007dagger}, based on examples from functional analysis, quantum mechanics, and  classical and quantum computing:
a \define{dagger category}  is
a category $\cC$  equipped with
 a bijection $(-)^\dagger \colon \hom(X,Y) \to \hom(Y,X)$ for each pair $(X,Y)$ of objects,
such that $(f^\dagger)^\dagger = f$ and $(f\circ g)^\dagger = g^\dagger \circ f^\dagger$ for all composable morphisms $f,g$.
In other words, $(-)^\dagger$ is a functor $\cC \to \cC^\op$ which is involutive and the identity on (the set of) objects.

Selinger's definition, by requiring the operation $\dagger$ to be the identity on the set of objects, is seemingly 
incompatible with a fundamental principle of (higher) category theory: one should never require two objects to be equal, but rather isomorphic. Category theorists hold this principle in such high regard that structures violating it are often referred to as \emph{evil}. Our definition of dagger higher category will be based on 
a formulation of ``dagger category'' which is equivalent, in the categorical sense, to Selinger's definition but which avoids the ``evil'' aspects of the definition.
Our definition is further supported and motivated by examples that have been observed in the study of $\rm C^*$ or $\rm W^*$-categories~\cite{2411.01678}, categorifications of the concept of Hilbert spaces~\cite{MR1448713, bartlett}, or unitary fusion categories~\cite{MR4133163}.
Moreover, using our definition, in this paper we construct a higher dagger structure on any extended bordism category (with stable tangential structure), allowing us to give the first complete definition of extended unitary topological quantum field theory.

Our discussion will focus on the essential ideas and will be short on technical details and essentially devoid of proofs; for this reason, we will refer to ``statements'' rather than ``lemmas'' or ``theorems,'' and strictly speaking all statements herein (other than those cited to the published literature) carry the ontological status of conjectures. We also list some true ``conjectures,'' of whose veracity we are less confident.  Our goal is to sketch the big picture, and we hope that our article can serve as a blueprint for a more thorough development of the subject. 

After the first version of this article appeared a year ago, our definitions have been applied and further developed in the study of unitary representation theory of generalized symmetries~\cite{Bartsch:2024ech,Bartsch:2025drc}, the description of topological defects and their gauging~\cite{2504.17764,2505.04761}, as well as in relation to higher Hilbert spaces~\cite{2410.05120}. 

\subsection*{Outline}
In Section~\ref{Sec: 1dagger} we recall the coherent formulation of dagger 1-categories and discuss a slight generalization to $(\infty,1)$-categories. In Section~\ref{sec: n dagger} we present our definition of dagger $n$-categories using the concept of flagged categories introduced by Ayala and Francis~\cite{MR3869643}, which we reformulate in terms of enriched category theory in Section~\ref{sec:enrichment}. In Section~\ref{sec.defbicat} we specialize to the case $n=2$. In Section~\ref{sec:unitaryduality} we discuss the interplay between dagger structures and the categorical notions of duality and adjoints. Section~\ref{sec:bord} outlines an application of higher dagger categories to quantum field theory, which was one of the motivations for the present work.

\section{Dagger \texorpdfstring{$(\infty,1)$}{(infty,1)}-categories}\label{Sec: 1dagger}

As mentioned in the introduction, the traditional definition of dagger categories is \define{evil} in the following sense \cite{daggerevil}: 
since  it explicitly references the set of objects in~$\cC$, dagger structures cannot be readily transported along general equivalences of categories. In detail: given a dagger category $(\cC, (-)^{\dagger})$ and an equivalence $F \colon \cC \simeq \cD$ with (weak) inverse $F^{-1}$, the transported ``dagger'' functor\footnote{A different formula for the transported structure supplies a dagger structure on $\cD$ such that $F^{-1}$ is a dagger functor~\cite[Lemma 2.1.16]{karvonen:dagger}, but then $F$ will typically not only fail to be a dagger equivalence, it will fail to be compatible with the dagger structures at all.}
$F \circ (-)^\dagger \circ F^{-1} \colon \cD \to \cD^\op$ is naturally part of a \define{weak involution} in the sense that its square is coherently-isomorphic to $\id_\cD$, but it will almost never be an identity on objects.
We violated the general rule of thumb: from a category, one can coherently extract its \emph{groupoid} of objects up to isomorphism, but not its set of objects. Non-evil notions refer only to coherently-extractable data.

Nevertheless, there is a well-developed ``dagger category theory'' that parallels the usual theory of categories, with dagger versions of functor, natural transformation, and the like~\cite{MR0808930,MR1444286,MR2091457,MR3663592,MR3687214,karvonen:dagger,MR3584697,MR3968429,srinivasan:logic}. 
From the perspective of dagger category theory, the coherently-extractable space of objects is the groupoid of objects and \emph{unitary} isomorphisms.

With this in mind, 
\cite{luukjan} developed the following ideas.
Consider the $(2,1)$-category $\Cat$ of $1$-categories, functors and natural isomorphisms. This $(2,1)$-category carries a $\bZ/2\bZ$-action which sends $\cC \mapsto \cC^\op$. An \emph{anti-involutive category} is a homotopy fixed point for this action: explicitly, it consists of a category $\cC$, a categorical equivalence $\dagger \colon \cC \overset\sim\to \cC^\op$, and a natural isomorphism $\eta \colon \dagger^\op \circ \dagger \overset\sim\Rightarrow \id_\cC$, such that for each object $x \in \cC$, $\eta_{x^\dagger}^{-1} = (\eta_x)^\dagger$.

If $\cG$ is a groupoid, then $\cG$ and $\cG^{\op}$ are canonically equivalent via the functor that is the identity on objects and sends morphisms to their inverses. 
In other words, the $\bZ/2\bZ$-action $\cC \mapsto \cC^\op$ on $\Cat$ trivializes when restricted to the full sub-$(2,1)$-category $\Gpd \subset \Cat$ of groupoids. 
In particular, given an anti-involutive category $\cC$, the groupoid $\iota_0 \cC$ of objects of $\cC$ inherits a coherent action by $\bZ/2\bZ$. 
Let $(\iota_0 \cC)^{\bZ/2\bZ}$ denote the groupoid of homotopy fixed points for this action. 

As a motivating example, consider $\cC = \cat{Vec}^\fd_\bC$  the category of finite-dimensional vector spaces, and $\dagger \colon V \mapsto \overline{V}^\vee$  the functor that
assigns the complex-conjugate  dual of a vector space. Then $(\iota_0 \cC)^{\bZ/2\bZ}$ is the groupoid of finite-dimensional \emph{Hermitian} vector spaces, i.e.\ vector spaces equipped with a nondegenerate conjugate-symmetric sesquilinear form which might not be positive-definite, and ``unitary'' isomorphisms thereof. In order to encode the theory of Hilbert spaces, we must specify extra data: we must mark some of the Hermitian spaces as preferred, 
and leave out the others. This marking  restricts the set of objects of $(\iota_0 \cC)^{\bZ/2\bZ}$, but keeps all unitary isomorphisms between them; in other words, the marked subspace $\cC_0 \subset (\iota_0 \cC)^{\bZ/2\bZ}$ should be a full subgroupoid. If there were a vector space that did not admit a Hilbert structure, then we should have already excised it from $\cC$; in other words, the composition $\cC_0 \subset (\iota_0 \cC)^{\bZ/2\bZ} \to \cC$ should be essentially surjective.

\begin{definition}\label{defn.dagger1cat}
  A \define{coherent dagger 1-category} is an anti-involutive category $(\cC,\dagger)$ together with a fully faithful subgroupoid $\cC_0 \hookrightarrow (\iota_0 \cC)^{\bZ/2\bZ}$ such that the induced map $\cC_0 \hookrightarrow (\iota_0 \cC)^{\bZ/2\bZ} \to \iota_0 \cC$ is essentially surjective.
 A morphism between coherent dagger 1-categories is a functor respecting the anti-involutions and the subgroupoids.
\end{definition}

Coherent dagger structures are non-evil: 
 if $\cC$ is a dagger category and $\cC \cong \cD$ an equivalence with an ordinary category, then $\cD$ acquires a natural coherent dagger structure.
The main theorem of \cite{luukjan} says that coherent dagger categories are equivalent in the appropriate sense to dagger categories as traditionally defined. Namely, coherent dagger categories can be \define{strictified} to dagger categories as traditionally defined:

\begin{theorem}[\cite{luukjan}]\label{luuk-strictification}
    Any dagger category $\cC$ defines a coherent dagger category by keeping $\cC$ and $\dagger$ as-is, and setting $\cC_0$ to be the groupoid of objects in $\cC$ and unitary morphisms between them. 
    The so-defined functor from the $(2,1)$-category of dagger categories, dagger functors, and unitary natural isomorphisms to the $(2,1)$-category of coherent dagger categories is an equivalence.
    
    Its inverse assigns to a coherent dagger category $(\cC, \dagger, \cC_0 \hookrightarrow (\iota_0 \cC)^{\bZ/2\bZ})$ the category whose objects are the objects of $\cC_0$ and whose morphisms are the morphisms in $\cC$ between the images of objects; this new category is equivalent to $\cC$ because of the essential surjectivity axiom, and its anti-involution is strictly trivial on the set of objects.
\end{theorem}

Definition~\ref{defn.dagger1cat} immediately generalizes to the $(\infty,1)$-world, where it was independently proposed by \cite{simonhenry}. An \define{$(\infty,1)$-category} is a homotopically-coherent version of ``category enriched in spaces.'' As is standard in higher category theory, we will interchangeably use the word \define{$\infty$-groupoid} and \define{space}
for ``nice topological space considered up to homotopy.'' Then an $(\infty,1)$-category has not a set of morphisms between any two objects, but rather a space of morphisms; composition is associative up to parameterized homotopy; this homotopy itself satisfies higher and higher homotopies relating different parenthesizations. There are many ways of axiomatizing the notion of ``$(\infty,1)$-category,'' and we refer the reader to~\cite{MR2664620,MR3466443}
for nice surveys.

Ignoring size issues, there is an $(\infty,1)$-category $\Cat_{(\infty,1)}$ whose objects are $(\infty,1)$-categories and whose morphisms are $(\infty,1)$-functors. There is a functor $\iota_0 \colon \Cat_{(\infty,1)} \to \cat{Space}$, corepresented by the terminal $(\infty,1)$-category $\{\pt\}$, which assigns to each $(\infty,1)$-category its space of objects. Moreover, $\Cat_{(\infty,1)}$ carries a canonical $\bZ/2\bZ$-action sending $\cC \mapsto \cC^\op$ \cite{MR2182378}.
An \emph{anti-involutive $(\infty,1)$-category} is a (homotopy) fixed point for this action. 

As in the 1-categorical case, the $\bZ/2\bZ$-action on $\Cat_{(\infty,1)}$ trivializes on the full subcategory $\cat{Space} \subseteq \Cat_{(\infty,1)}$. Hence, the space of objects $\iota_0 \cC$ of any anti-involutive $(\infty,1)$-category $\cC$ inherits a coherent $\bZ/2\bZ$ action, whose homotopy fixed points we continue to denote by $(\iota_0 \cC)^{\bZ/2\bZ}. $

\begin{definition}\label{defn.dagger1catredux}
A \define{dagger $(\infty,1)$-category} is an anti-involutive $(\infty,1)$-category, together with a full sub-$\infty$-groupoid $\cC_0 \hookrightarrow (\iota_0 \cC)^{\bZ/2\bZ}$ such that the induced map $\cC_0 \hookrightarrow (\iota_0 \cC)^{\bZ/2\bZ} \to \cC$ is essentially surjective.
\end{definition}

Anticipating our $n$-dimensional generalization, let us note that part of the data of a dagger $(\infty,1)$-category is the space $\cC_0$, the $(\infty,1)$-category $\cC$, and the essential surjection $\cC_0 \to \iota_0 \cC$. This structure is called a \define{flagged $(\infty,1)$-category}. Any specific $1$-category $\cC$ in the traditional sense supplies a flagged $(\infty,1)$-category: one takes $\cC_0$ to be the actual set of objects, whereas $\iota_0 \cC$ is the groupoid of objects and isomorphisms. In other words, flaggings are a way of remembering  ``sets'' of objects in a way that nevertheless transports coherently. The requirement that $\cC_0 \to \iota_0 \cC$ be essentially surjective simply means that the data of the homomorphisms between elements of $\cC_0$ suffices to recover all of $\cC$ up to categorical equivalence.
 A flagging is called \define{univalent}\footnote{This is reminiscent of the completeness condition of Segal spaces.} if it does not in fact remember any further data: if the map $\cC_0 \to \iota_0 \cC$ is not just essentially surjective but also fully faithful.
 
 By analogy, in a dagger $(\infty,1)$-category, the requirement that $\cC_0 \to (\iota_0 \cC)^{\bZ/2\bZ}$ be fully faithful will be called the \define{univalence axiom}, even though it does not force any map to be an equivalence. And by analogy, dropping the univalence axiom leads to a useful weakening of Definition~\ref{defn.dagger1catredux}. 
 \begin{definition}
  A \define{flagged dagger $(\infty,1)$-category} is an anti-involutive $(\infty,1)$-category $\cC$ equipped with a map of spaces $\cC_0 \to (\iota_0 \cC)^{\bZ/2\bZ}$ such that the induced map $\cC_0 \to (\iota_0 \cC)^{\bZ/2\bZ} \to \cC$ is essentially surjective.
 \end{definition}

 Our motivating example of a flagged dagger $(\infty,1)$-category is the higher category of bordisms, which we discuss in Example~\ref{eg:bordism}.
 
Whether univalent or not, the flagging should be thought of as recording the ``identity on objects'' condition in the traditional definition of dagger category.
In a flagged dagger category which is not univalent, the groupoid $\cC_0$ selects a notion of equivalence between objects that is finer than unitary equivalence.

Write $\flagdagCat_{(\infty,1)}$ for the $(\infty,1)$-category of flagged dagger $(\infty,1)$-categories: the homomorphisms are the obvious ones which preserve the equivariance and the flaggings. Write $\dagCat_{(\infty,1)}$ for the full subcategory of $\cat{Fl}\dagCat_{(\infty,1)}$ on the univalent ones. 

\begin{statement}\label{statement:completionflagging}
The inclusion
$$ \dagCat_{(\infty,1)} \hookrightarrow \flagdagCat_{(\infty,1)}$$
 admits a left adjoint, which replaces $\cC_0$ with the full subgroupoid inside $(\iota_0 \cC)^{\bZ/2\bZ}$ on the essential image of $\cC_0 \to (\iota_0 \cC)^{\bZ/2\bZ}$.
 \end{statement}

One benefit of Statement~\ref{statement:completionflagging} is that it can be easier to present flagged dagger categories than their univalentizations. 
Indeed, most $(\infty,1)$-categories $\cC$ of interest
(for example the bordism categories discussed in Section~\ref{sec:bord})
 come with distinguished flaggings $\cC_0 \to \cC$, and sometimes this flagged category is naturally anti-involutive in the sense that $\cC$ is anti-involutive and $\cC_0 \to \iota_0 \cC$ is $\bZ/2\bZ$-equivariant. To promote this to a flagged dagger structure then merely requires a trivialization of the $\bZ/2\bZ$-action on $\cC_0$. Speaking very approximately, it can be easier to do this when $\cC_0$ has very few morphisms requiring trivialization data.

Instead of dropping the univalence axiom, we could instead drop the essential surjectivity condition:
\begin{definition}\label{defn:coflagged}
  A \define{coflagged dagger $(\infty,1)$-category} is an anti-involutive $(\infty,1)$-category $\cC$ equipped with a fully faithful inclusion of spaces $\cC_0 \hookrightarrow (\iota_0 \cC)^{\bZ/2\bZ}$.
\end{definition}

Comparing to Definition~\ref{defn.dagger1cat}, a coflagged dagger $(\infty,1)$-category can be thought of as an anti-involutive $(\infty,1)$-category equipped with the data of Hermitian pairings on some but not necessarily all objects.

\begin{statement}\label{statement:completioncoflagging}
  Writing $\coflagdagCat_{(\infty,1)}$ for the $(\infty,1)$-category of coflagged dagger $(\infty,1)$-categories, the inclusion
  $$ \dagCat_{(\infty,1)} \hookrightarrow \coflagdagCat_{(\infty,1)} $$
  admits a right adjoint, which replaces $\cC$ with its full subcategory on the image of $\cC_0 \to (\iota_0 \cC)^{\bZ/2\bZ} \to \cC$.
\end{statement}
Note that every anti-involutive $(\infty,1)$-category is naturally a coflagged dagger $(\infty,1)$-category with $\cC_0 = (\iota_0 \cC)^{\bZ/2\bZ}$.
Through this observation, Statement~\ref{statement:completioncoflagging} specializes to a higher version of the ``Hermitian completion'' functor of \cite{luukjan}.

\section{Dagger \texorpdfstring{$n$}{n}-categories via flaggings}\label{sec: n dagger}

The notion of $(\infty,n)$-category is designed to formalize categories with an $n$-fold hierarchy of directions of composition. As with the $(\infty,1)$-case, there are many models, some of which are surveyed in \cite{MR3109865,MR4186138,MR4301559}. Informally, an $(\infty,n)$-category has a space of objects; for each pair of objects, a space of $1$-morphisms between them and a coherently-associative and coherently-unital composition law; for each pair of parallel $1$-morphisms, a space of $2$-morphisms between them, and two (coherently-associative and coherently-unital) composition laws, one of which lifts the composition law of $1$-morphisms and the other of which is in the new second dimension; and so on up to dimension $n$. 

In this informal description of $(\infty,n)$-categories in the previous paragraph, we did not specify which maps between them should be considered equivalences. There is a most natural guess if one explained the idea without supplying the name ``category'': a homomorphism could be declared an equivalence when it induces equivalences on all spaces of $k$-morphisms, including on the space of objects. Comparing with the traditional notion of strict 1-category, this choice would select the strict isomorphisms and not the categorical equivalences. In contrast, categorical equivalence does not remember any specific sets or spaces of $k$-morphisms, but merely the higher groupoids thereof. As in the $(\infty,1)$-case, the extra data of  spaces of morphisms in a presentation of an $(\infty,n)$-category is called a \define{flagging} of that $(\infty,n)$-category in \cite{MR3869643}, where such flaggings are introduced as model-independent interpretations of Segal sheaves on Joyal's category $\Theta_n$ (or equivalently of $n$-fold Segal spaces). Given a correct theory of $(\infty,n)$-categories, the notion of flagged $(\infty,n)$-category can be defined as follows:

\begin{definition}[{\cite[Def.~0.12]{MR3869643}}]
  A \define{flagged $(\infty,n)$-category} is a chain
  $$ \cC_0 \to \cC_1 \to \dots \to \cC_n$$
  where each $\cC_k$ is an $(\infty,k)$-category, and the map $\cC_k \to \cC_{k+1}$ is essentially surjective on $(\leq k)$-morphisms\footnote{A functor $F \colon \cC \to \cD$ of $(\infty,n)$-categories is \define{essentially surjective on $(\leq k)$-morphisms} if it is essentially surjective on objects and, 
  for every $0 \leq j < k$ and every pair of $j$-morphisms $f$ and $g$ with the same source and target, the functor $F_{f,g} \colon \hom_\cC(f,g) \to \hom_\cD(Ff, Fg)$ of $(\infty,n-j-1)$-categories is essentially surjective.}. The flagging is \define{univalent} if $\cC_k \to \iota_k \cC_{k+1}$ is an equivalence for each $k$, where $\iota_k \colon \Cat_{(\infty,k+1)} \to \Cat_{(\infty,k)}$ is the right-adjoint to the inclusion $\Cat_{(\infty,k)} \hookrightarrow \Cat_{(\infty,k+1)}$, which takes an $(\infty,k+1)$-category and forms an $(\infty,k)$-category by forgetting the non-invertible $(k+1)$-morphisms.
\end{definition}

The definition of dagger $(\infty,1)$-category used as an essential ingredient the fact that every $(\infty,1)$-category $\cC$ has an \define{opposite} $\cC^\op$ that reverses the direction of composition. Similarly, 
given an $(\infty,n)$-category $\cC$, one can produce new $(\infty,n)$-categories by reversing any of the $n$ directions of composition. This supplies an action of $(\bZ/2\bZ)^n$ on $\Cat_{(\infty,n)}$, the $(\infty,1)$-category of $(\infty,n)$-categories. In fact, this action is completely canonical: one of the theorems of \cite{MR4301559}, generalizing \cite{MR2182378}, says that this map
\begin{equation} \label{eqn:unicity}
    (\bZ/2\bZ)^n \to \Aut(\Cat_{(\infty,n)})
\end{equation} 
is an isomorphism of (higher) groups. Given $1 \leq k \leq n$, we will occasionally write $(\bZ/2\bZ)_k \subset (\bZ/2\bZ)^n$ for the $k$th coordinate $\bZ/2\bZ$, which acts on $(\infty,n)$-categories by reversing the composition of $k$-morphisms; we will write that action as $\cC \mapsto \cC^{k\op}$.

\begin{definition}\label{defn:volutive}
 Given a group homomorphism $G \to (\bZ/2\bZ)^n$, an $(\infty,n)$-category $\cC  \in \Cat_{(\infty,n)}$ is 
 \define{$G$-volutive}\footnote{From the Latin \emph{volvere}, ``to roll,'' and \emph{involvere}, ``to roll inwards.''} 
 when it is equipped with the data making it into a fixed point for the action of $G$ on $\Aut(\Cat_{(\infty,n)})$ via~\eqref{eqn:unicity}. 
    In the special case when $G \to (\bZ/2\bZ)^n$ is the identity, we will say $\cC$ is \define{fully-volutive}. When $G = (\bZ/2\bZ) \to (\bZ/2\bZ)^n$ selects the last involution, we will say $\cC$ is \define{top-\wipedge}.
\end{definition}

As in the 1-category case, \wipedge\ structures do not capture the theory of daggers: they do not include an analogue of the requirement that ``dagger is the identity on objects.'' To encode the latter, we use flaggings. To set up the definition, we note the following. Every $(\infty,k)$-category is in particular an $(\infty,k+1)$-category, and the inclusion $\Cat_{(\infty,k)} \hookrightarrow \Cat_{(\infty,k+1)}$ is stable under the ambient $(\bZ/2\bZ)^{k+1}$-action. Indeed, the first $k$ involutions $(\bZ/2\bZ)^k \subset (\bZ/2\bZ)^{k+1}$ act on $\Cat_{(\infty,k)}$ via the canonical action, and the last involution $(\bZ/2\bZ)_{k+1}$ has a canonical trivialization.
Recall that a fixed point for the trivial action of a group $G$ on a category $\cX$ is precisely an action of $G$ on some object $X \in \cX$. Thus if $\cC$ is an $(\infty,k)$-category thought of as an $(\infty,k+1)$-category, then a fully-\wipedge\ structure in the $(\infty,k+1)$-sense is precisely a fully-\wipedge\ structure in the $(\infty,k)$-sense together with a coherently-compatible $\bZ/2\bZ$-action. In particular, by picking the trivial $\bZ/2\bZ$-action, any fully-\wipedge\ $(\infty,k)$-category can be thought of as a fully-\wipedge\ $(\infty,k+1)$-category.

\begin{definition}\label{defn.flagdagncat}
  A \define{flagged fully-dagger $(\infty,n)$-category} is a flagged $(\infty,n)$-category
  $$ \cC_0 \to \cC_1 \to \dots \to \cC_n$$
  such that each $\cC_k$ is a fully-\wipedge\ $(\infty,k)$-category, and the map $\cC_k \to \cC_{k+1}$ is a map of fully-\wipedge\ $(\infty,k+1)$-categories, where $\cC_k$ is given the trivial $(k+1)$th anti-involution. Equivalently, $\cC_k \to \cC_{k+1}$ is factored through a map $\cC_k \to (\iota_k \cC_{k+1})^{(\bZ/2\bZ)_{k+1}}$ of fully-\wipedge\ $(\infty,k)$-categories.
\end{definition}

We remind the reader that, in the higher-categorical world, requesting that ``this is a map of these things'' is requesting for extra structure on the map. In the case at hand, this structure can be unpacked as follows. A map $\cC_k \to \cC_{k+1}$ is the same as a map $\cC_k \to \iota_k \cC_{k+1}$. The fully-\wipedge\ structure on $\cC_{k+1}$ induces a fully-\wipedge\ structure on $\iota_k \cC_{k+1}$ in the $(\infty,k)$-sense together with a typically-nontrivial action of $\bZ/2\bZ$. The requested structure unpacks to a map of fully-\wipedge\ $(\infty,k)$-categories $\cC_k \to (\iota_k \cC_{k+1})^{\bZ/2\bZ}$. Comparing with Definitions~\ref{defn.dagger1cat} and~\ref{defn.dagger1catredux}, we exactly recover the ``flagged dagger categories''; in particular, the ``essential surjectivity'' axiom is enforced by asking that $ \cC_0 \to \cC_1 \to \dots \to \cC_n$ be a flagged $(\infty,n)$-category independent of the \wipedge, but we do not have any univalence axiom. We add that univalence axiom now.

\begin{definition}\label{defn.dagncat}
  A flagged fully-dagger $(\infty,n)$-category is \define{univalent} if the maps $\cC_k \to (\iota_k \cC_{k+1})^{\bZ/2\bZ}$ are fully-faithful on $(\geq k+1)$-morphisms\footnote{A functor $F\colon \cC \to \cD $ is \define{fully-faithful on $(\geq k+1)$-morphisms} if, for every pair of $k$-morphisms $f$ and $g$ with the same source and target, the functor $F_{f,g}\colon \cC(f,g)\to \cD(F(f),F(g))$ is an equivalence.}.
  A \define{fully-dagger $(\infty,n)$-category} is a univalent flagged fully-dagger $(\infty,n)$-category.  We will write $\dagCat_{(\infty,n)}$ for the $(\infty,1)$-category of fully-dagger $(\infty,n)$-categories.
\end{definition}
Beware that a fully-dagger $(\infty,n)$-category is \emph{not} the same as a flagged fully-dagger $(\infty,n)$-category of which the underlying flagging is univalent.

Definitions~\ref{defn.flagdagncat} and~\ref{defn.dagncat} refer only to the successive inclusions $\cC_k \to (\iota_k \cC_{k+1})^{\bZ/2\bZ} \to \cC_{k+1}$. But they imply more general conditions, and these more general conditions will be needed when studying groups $G$ other than $(\bZ/2\bZ)^n$.

\begin{statement}\label{statement:moregeneralfactorization}
  In a flagged fully-dagger
  $(\infty,n)$-category $\cC$, the map $\cC_k \to \cC_l$ is essentially surjective on $(\leq k)$-morphisms for every $0 \leq k \leq l \leq n$. This map factors through $(\iota_k \cC_{l})^{(\bZ/2\bZ)_{k+1}\times\dots\times (\bZ/2\bZ)_l}$. If $\cC$ is univalent, then $\cC_k \to (\iota_k \cC_{l})^{(\bZ/2\bZ)_{k+1}\times\dots\times (\bZ/2\bZ)_l}$ is fully faithful on $(\geq l)$-morphisms.
\end{statement}

\begin{definition}\label{defn:G-dagger-nonadjunct}
  Let $G \subset (\bZ/2\bZ)^{n}$ be a subgroup. For any subinterval $\{k+1,\dots,l\} \subset \{1,\dots,n\}$, define $G(\{k+1,\dots,l\}) = G \cap_{(\bZ/2\bZ)^n} \prod_{j=k+1}^l (\bZ/2\bZ)_{j}$. A \define{flagged $G$-dagger $(\infty,n)$-category} is a flagged $(\infty,n)$-category $\cC$ together with, for each $k$, a $G(\{1,\dots,k\})$-\wipedge\ structure on $\cC_k$ and, for each $0 \leq k< l \leq n$, and a factorization of $\cC_k \to \cC_l$ through a map $\cC_k \to (\iota_k \cC_l)^{G(\{k+1,\dots,l\})}$ of $G(\{1,\dots,k\})$-\wipedge\ $(\infty,k)$-categories. A flagged $G$-dagger $(\infty,n)$-category is \define{univalent} when the maps $\cC_k \to (\iota_k \cC_l)^{G(\{k+1,\dots,l\})}$ are all fully faithful on $(\geq l)$-morphisms, in which case the flagged $G$-dagger $(\infty,n)$-category is a \define{$G$-dagger $(\infty,n)$-category}.
\end{definition}

We will write $G\dagCat_{(\infty,n)}$ for the $(\infty,1)$-category of $G$-dagger $(\infty,n)$-categories.
The most interesting case for examples is when $G = (\bZ/2\bZ)_n$, in which case we will also refer to $(\bZ/2\bZ)_n$-dagger $(\infty,n)$-categories as \define{top-dagger}. More generally, the cases that arise in examples are when $G \subset (\bZ/2\bZ)^n$ is a product of some of the coordinate $(\bZ/2\bZ)_k$s. For such a $G$, Definition~\ref{defn:G-dagger-nonadjunct} simplifies: since in this case the  maps $\prod_{j={k+1}}^l G(\{j\}) \to G(\{k+1,\dots,l\})$ are all isomorphisms,
one can restrict to just the successive inclusions $\cC_k \to \cC_{k+1}$ in the definition and invoke a $G$-version of Statement~\ref{statement:moregeneralfactorization}. We have suggested a more general definition in anticipation of our discussion of unitary duality in Section~\ref{sec:unitaryduality}.

\section{Dagger \texorpdfstring{$n$}{n}-categories via enrichment}
\label{sec:enrichment}

Let $\cV$ be a symmetric monoidal $(\infty,1)$-category, with monoidal structure $\otimes_\cV \colon \cV \times \cV \to \cV$ and monoidal unit $1_\cV \in \cV$. There is a notion of \define{$\cV$-enriched $(\infty,1)$-category}, which we will abbreviate as \define{$\cV$-category}. A full definition is provided by \cite{MR3345192}; we will recall the main ingredients. A \define{flagged $\cV$-category} $\cC$ consists of a space $\cC_0$ of objects and, for each $x,y \in \cC_0$, an object $\hom(x,y) \in \cV$, together with a coherently-associative and coherently-unital composition law $\hom(y,z) \otimes_\cV \hom(x,y) \to \hom(x,z)$. The \define{global sections} functor $\Gamma := \hom(1_\cV, -) \colon \cV \to \cat{Space}$ is lax-symmetric-monoidal, and so $\cC$ induces a flagged plain $(\infty,1)$-category $\Gamma\cC$ with the same space $\cC_0$ of objects, and with spaces of morphisms given by $\Gamma \hom(-,-)$. The flagged $\cV$-category $\cC$ is called \define{univalent} if $\Gamma\cC$ is univalent. By definition, a \define{$\cV$-category} is a univalent flagged $\cV$-category. We will write $\Cat[\cV]$ for the $(\infty,1)$-category of $\cV$-categories. The primordial example: $\Cat[\Cat_{(\infty,n-1)}] \cong \Cat_{(\infty,n)}$.

Using the symmetry on $\cV$, for each $\cV$-category $\cC$ it should be possible to define an \define{opposite} $\cV$-category $\cC^\op$ with the same objects and morphisms but the opposite order of composition. This supplies an involution $\bZ/2\bZ \to \Aut(\Cat[\cV])$. We can immediately generalize Definition~\ref{defn.dagger1catredux}:
\begin{definition}\label{defn.enricheddagger}
  Let $\cC$ be a $\cV$-category, with space of objects $\iota_0 \cC$. A \define{flagged dagger structure} on $\cC$ is a fixed-point structure on $\cC$ for the $\bZ/2\bZ$-action on $\Cat[\cV]$ together with a map of spaces $\cC_0 \to (\iota_0\cC)^{\bZ/2\bZ}$ such that the induced map $\cC_0 \to (\iota_0\cC)^{\bZ/2\bZ} \to \iota_0 \cC$ is essentially surjective. A flagged dagger structure is \define{univalent} when $\cC_0 \to (\iota_0\cC)^{\bZ/2\bZ}$ is fully faithful. A \define{dagger $\cV$-category} is a $\cV$-category with a univalent flagged dagger structure.
\end{definition}
We will write $\dagCat[\cV]$ for the $(\infty,1)$-category of dagger $\cV$-categories. We expect that the
equivalence $\Cat[\Cat_{(\infty,n-1)}] \cong \Cat_{(\infty,n)}$  generalizes to dagger categories:

\begin{statement}\label{thm.comparison}
  There is an  equivalence of $(\infty,1)$-categories
  $$ \dagCat[\dagCat_{(\infty,n-1)}] \cong \dagCat_{(\infty,n)}.$$
\end{statement}

In particular, by iterating Statement~\ref{thm.comparison}, one arrives at an alternative model of  {fully-dagger $(\infty,n)$-category} than the one given in Definition~\ref{defn.dagncat}. Mixing $\Cat[-]$ and $\dagCat[-]$ provides alternative models of the versions with dagger structures on only some levels:

\begin{statement}
\label{thm.enrichedtopdagger}
  There are equivalences of $(\infty,1)$-categories
  \begin{gather*}
    \dagCat[\Cat_{(\infty,n-1)}] \cong (\bZ/2\bZ)_1\dagCat_{(\infty,n)} \\
    \Cat[G\dagCat_{(\infty,n-1)}] \cong (* \times G)\dagCat_{(\infty,n)}
  \end{gather*}
\end{statement}
For example, iterating the second of these equivalences supplies an alternative model for top-dagger $(\infty,n)$-categories.

\section{Dagger bicategories}\label{sec.defbicat}

In this section, we will elaborate on how Definitions~\ref{defn.dagncat} and~\ref{defn:G-dagger-nonadjunct} play out in the case of bicategories. 
Our goal is to outline a strictification  for coherent dagger bicategories analogous to Statement~\ref{luuk-strictification} for $1$-categories.
Even though our definition of a fully-dagger bicategory unpacks to something  complicated, 
it can be strictified to a  ``more traditionally defined'' fully-dagger bicategory involving less data (but more ``evil''):

\begin{definition}
\label{def:biinvbicat}
    A \define{bi-involutive bicategory} is a bicategory $\mathcal{B}$ equipped with two functors
    \begin{align*}
        \dagger_1\colon \cB \to \cB^{1\op}, \quad \dagger_2\colon \cB \to \cB^{2\op}
    \end{align*}
    such that 
    \begin{enumerate}
        \item $\dagger_2$ is the identity on objects and $1$-morphisms and strictly squares to the identity;
        \item $\dagger_1$ is the identity on objects (but not necessarily on $1$-morphisms) and weakly squares to the identity in the sense that it comes equipped with a natural isomorphism $\phi\colon \dagger_1^2 \to \id_{\mathcal{B}}$, which is the identity on objects (but not necessarily on $1$-morphisms).
    \end{enumerate}
    There is a further condition, left to the reader, comparing the two ways to trivialize $\dagger_1^3$.
    Additionally, the two daggers strongly commute:
    \[
    \dagger_2^{1\op}\circ \dagger_1 = \dagger_1^{2\op} \circ \dagger_2.
    \]
    This equality should be compatible with the isomorphism $\phi \colon \dagger_1^2 \to \id_{\mathcal{B}}$, i.e.\ $\phi$ is unitary with respect to $\dagger_2$. 
\end{definition}

Unpacking Definition~\ref{def:biinvbicat}, it is an enriched-type definition in the sense of Section~\ref{sec:enrichment}, so that Statement~\ref{thm.comparison} for $n=2$ gives a justification for the hope that bi-involutive bicategories are a model for fully dagger bicategories.
In other words, similarly to how a bicategory is a category weakly enriched in categories, a bi-involutive bicategory is a $\dagger$-category ($\dagger_1$) weakly enriched in $\dagger$-categories ($\dagger_2$).

We decided on the name ``bi-involutive bicategory'' because it generalizes the notion of a bi-involutive tensor category of \cite{MR3663592} to a bicategory with more than one object. 

\begin{remark}\label{rem:top and bottom daggers}
There are at least two other structures one might call $\dagger$-bicategories, which are both $G$-dagger categories in the sense of~\ref{defn:G-dagger-nonadjunct}:
they correspond to considering the two canonical $\bZ/2\bZ$-subgroups of $\bZ/2\bZ\times \bZ/2\bZ$.\footnote{It could be interesting to consider the diagonal subgroup.} 
The first one leads to a $\dagger_2$ which is the identity on objects and 1-morphisms as above, giving a category weakly enriched in dagger categories as a special case of Statement~\ref{thm.enrichedtopdagger}.
The other structure leads to a $\dagger_1$ which is the identity on objects as above, a dagger category weakly enriched in categories. 
\end{remark}

\begin{example}
A natural example of a bi-involutive bicategory is the bicategory of von Neumann algebras, Hilbert space bimodules equipped with commuting normal actions
and bounded bimodule homomorphisms \cite{MR0703809,MR0945550,MR1303779,MR3342166,MR4419534}. 
The daggers are given by
\begin{itemize}
    \item[$\dag_1$] For a Hilbert space bimodule ${}_MH_N$ for $M,N$ von Neumann algebras, $H^{\dag_1}$ is defined to be ${}_N\overline H_M$ where $\overline H$ is the complex conjugate Hilbert space and the actions are given by $b\cdot\overline\xi\cdot a=\overline{a^*\xi b^*}$
    \item[$\dag_2$] For 2-morphisms (Hilbert space bimodule homomorphisms), the dagger is simply given by the adjoint as maps between Hilbert spaces.
\end{itemize}
We note that there is actually an involution on objects (which we may call $\dag_0$) given by taking the opposite algebra (or equivalently, the complex conjugate algebra). 
One reason for this additional object-level involution in monoidal bicategories is due to the extra $\mathbb Z/2\bZ$ action given by reversing the monoidal product.
\end{example}

\begin{example}
    Consider the case in which $\cB$ has one object so that $\cB$ is the delooping of a monoidal category $\cC$.
    Having only a top-dagger $\dagger_2$ on $\cB$ is the same as having a monoidal dagger structure on $\cC$.
On the other hand, having only a $\dagger_1$ on $\cB$ is the same as having a weak covariant involution $\overline{(.)}$ on $\cC$ that is op-monoidal in the sense that
\(
\overline{x \otimes y} \cong \overline{y} \otimes \overline{x}.
\)
Such structures have been considered in the context of dagger categories in~\cite{MR2861112}.
\end{example}

In the rest of this section, we provide some ideas for a proof of the following strictification result, which is a bicategorical analogue of the main theorem of \cite{luukjan}.

\begin{statement}
    The $(3,1)$-category of bi-involutive bicategories is equivalent to the $(3,1)$-category of fully-dagger bicategories.
\end{statement}

Firstly, a fully-\wipedge\ structure on a bicategory $\mathcal{B}$ consists of a pair of equivalences $\psi_1\colon \mathcal{B}\to \mathcal{B}^{1\op}$ and $\psi_2\colon \mathcal{B}\to \mathcal{B}^{2\op}$, together with natural isomorphisms $ \Omega_1\colon \psi_1^{1\op} \circ \psi_1 \to \id_{\mathcal{B}}$, $ \Omega_2 \colon \psi_2^{2\op} \circ \psi_2 \to \id_{\mathcal{B}}$, and $\Omega_{12}\colon \psi_2^{1\op}\circ \psi_1 \to \psi_1^{2\op} \circ \psi_2$, and various modifications which satisfy various coherence conditions. 
To provide a fully-dagger structure on this \wipedge\ bicategory, we need to understand the relevant fixed point categories and specify our flaggings.

The bigroupoid $(\iota_0 \mathcal{B})^{\bZ/2\bZ\times \bZ/2\bZ}$ has objects consisting of an element of $\mathcal{B}$ equipped with fixed point data for both  $\psi_1$ and $\psi_2$ together with compatibility data between them, satisfying various coherence conditions.
Similarly to the $1$-categorical case, this in particular consists of isomorphisms $h_b^1\colon \psi_1(b) \to b$ and $h_b^2\colon \psi_2(b) \to b$.
The space $\cB_0$ in the definition of a dagger bicategory can be identified with a subspace of $(\iota_0 \mathcal{B})^{\bZ/2\bZ\times \bZ/2\bZ}$ such that the map to $\cB$ is essentially surjective. This picks out at least one preferred $\bZ/2\bZ\times \bZ/2\bZ$-fixed point structure for every object $b\in \mathcal{B}$. 

The bicategory $\iota_1\mathcal{B}^{\bZ/2\bZ}$ has objects given by an element of $\mathcal{B}$ equipped with fixed point data for $\psi_2$. The 1-morphisms consist of pairs of a 1-morphism $f\colon b\to b'$ in $\mathcal{B}$ and 2-isomorphisms 
\begin{equation}
\begin{tikzcd}
     \psi_2(b) \ar[r,"h_b^2"]  \ar[d, "\psi_2(f)",swap] & b \ar[d, "f"] \ar[ld, Rightarrow] \\ 
     \psi_2(b') \ar[r,"h_{b'}^2", swap] & b'
\end{tikzcd}
\label{Eq: hf}
\end{equation}
satisfying a natural coherence condition. 
One can think of this $2$-morphism as the data specifying how $f$ is a unitary $1$-morphism.
The bicategory $ \cB_1$ is equivalent to a subcategory of $\iota_1\mathcal{B}^{\bZ/2\bZ}$, whose objects can be identified with those of $\cB_0$, since $\cB_0\to \cB_1$ is essentially surjective. The 1-morphisms pick out at least one preferred $\bZ/2\bZ$-fixed point structure on every 1-morphism. 
The 2-morphisms are fixed by the condition that the map to $\iota_1\cB^{\bZ/2\bZ}$ is fully faithful on 2-morphisms. 
Note how this additionally fixes the 1-morphisms of $\cB_0$: they are exactly those whose fixed point structure restricts to one of the chosen fixed point structures in $\cB_1$. 
In summary, the structure of a dagger bicategory is equivalent to a fully-\wipedge\ bicategory $\mathcal{B}$ together with a choice of at least one $\bZ/2\bZ\times \bZ/2\bZ$-fixed point on every object $b\in \mathcal{B}$ and at least one compatible $\bZ/2\bZ$-fixed point on every 1-morphism. 

From this data we can construct a bi-involutive bicategory (analogous to the construction in \cite{luukjan}) as follows.
Consider the bicategory $\mathcal{B}'$ of which the objects are the objects of $\cB_0$ and 1-morphisms the 1-morphisms of $\cB_1$. 
The 2-morphisms are those of $\mathcal{B}$, requiring no compatibility with the fixed point data. 
By construction, there is a forgetful functor $\mathcal{B}'\to \mathcal{B}$, which is an equivalence of bicategories. 
The anti-involutions $\psi_1$ and $\psi_2$ induce compatible anti-involutions on $\mathcal{B}'$.\footnote{We expect the equivalence $\cB' \to \cB$ to come equipped with a canonical datum saying it preserves the fully-\wipedge\ structures.}
Namely, $\dagger_1$ is defined by the same formula as in the $1$-categorical case:
\begin{align*}
    f\colon b \to b' & \longmapsto h_b^1 \circ \psi_1(f) \circ (h_{b'}^1)^{-1}.
\end{align*} 
We do not spell out the fixed point structure \eqref{Eq: hf} on the $1$-morphism $f^{\dagger_1}$ here.
There is a natural isomorphism $\phi\colon \dagger_1^2\to \id_{\mathcal{B}'}$ which is the identity on objects and uses the fixed point data on $b$ and $b'$ to get a $2$-isomorphism $f^{\dagger_1 \dagger_1} \cong f$ for a $1$-morphism $f\colon b \to b'$.
The top-dagger $\dagger_2 \colon \mathcal{B}'\longrightarrow \mathcal{B}'^{2\op} $ is the identity on objects and 1-morphisms. It sends a 2-morphism $\vartheta\colon f \Rightarrow g$ to a version of $\psi_2(\vartheta)$, where we identify its domain and target with $g$ and $f$ respectively using their fixed point data.
The functor $\dagger_2$ strictly squares to the identity. 
The two daggers commute strictly, 
and $\phi$ is $\dagger_2$-unitary. 
We have thus constructed a canonical bi-involutive bicategory $\cB'$ from the coherent full $\dagger$ bicategory $\cB$.

\section{Dagger \texorpdfstring{$n$}{n}-categories with unitary duals}\label{sec:unitaryduality}

An $(\infty,n)$-category $\cC$ is said to
 \define{have adjoints}, if for each $1 \leq k \leq n-1$ and each $k$-morphism $f \colon x \to y$ (between parallel $(k-1)$-morphisms $x,y$), there exist $k$-morphisms $f^R,f^L \colon y \to x$ and $(k+1)$-morphisms $\eta^R\colon \id_x \to f^R\circ f$, $\epsilon^R \colon f \circ f^R \to \id_y$, $\eta^L\colon \id_y \to f\circ f^L$, $\epsilon^L \colon f^L \circ f \to \id_x$ such that the \define{zig-zag} compositions, which after suppressing coherence (e.g.\ unitor and associator) information become
 \begin{align*}
 f \overset{f\eta^R}\longrightarrow ff^R f \overset{\epsilon^R f}\longrightarrow f,&&
 f^R \overset{\eta^R f^R}\longrightarrow f^R f f^R \overset{f^R \epsilon^R}\longrightarrow  f^R, \\
 f \overset{\eta^L f}\longrightarrow ff^L f \overset{f\epsilon^L}\longrightarrow f,&&
 f^L \overset{f^L\eta^L}\longrightarrow f^L f f^L \overset{\epsilon^L f^L}\longrightarrow f^L,
\end{align*}
  are equivalent to identities. Let $\AdjCat_{(\infty,n)} \subset \Cat_{(\infty,n)}$ denote the full sub-$(\infty,1)$-category on the $(\infty,n)$-categories with adjoints.

The theory of $(\infty,n)$-categories with adjoints has been well-studied, but there are many questions remaining. The most famous work pertains to the symmetric monoidal case, in which case dualizability is also imposed on objects. We will call a symmetric monoidal category \define{rigid} if it has adjoints and also admits duals for objects. The \define{Cobordism Hypothesis} of \cite{MR1355899, MR2555928} asserts that the free rigid symmetric monoidal $(\infty,n)$-category generated by a single object is the $(\infty,n)$-category $\Bord_n^\fr$ of framed $n$-dimensional bordisms---a \define{framing} $\tau$ on a smooth $n$-manifold $M$ is a trivialization of its tangent bundle $\tau \colon \rT_M \cong \bR^n$. In other words, if $\cC$ is a rigid symmetric-monoidal $(\infty,n)$-category, with space of objects $\iota_0 \cC$, then there is a canonical equivalence
$$ \iota_0 \cC \cong \hom_{\mathrm{sym \otimes}}(\Bord_n^\fr, \cC).$$
The group $\rO(n)$ obviously acts on the framings on a given $n$-manifold, by rotating the trivialization, and so acts on the $(\infty,n)$-category $\Bord_n^\fr$. This in turn supplies a famous action of $\rO(n)$ on $\iota_0 \cC$ for any rigid symmetric monoidal $(\infty,n)$-category $\cC$.

In fact, there is a larger group that acts, as explained in \cite[Remark 2.4.30]{MR2555928}. It makes sense to talk about a ``tangent bundle'' of a piecewise linear (a.k.a. PL) manifold, but it is not a vector bundle: whereas the tangent bundle of a smooth $n$-manifold $M$ is classified by a map $\rT_M \colon M \to \rB \rO(n) \simeq \rB \Diff(n)$, a tangent bundle of a PL manifold is classified by a map $\rT_M \colon M \to \rB\PL(n)$. In particular, it makes sense to talk about \define{framed PL manifolds}: they are PL manifolds equipped with a trivialization of $\rT_M \colon M \to \rB\PL(n)$. By construction, any smoothing of a PL manifold lifts $\rT_M \colon M \to \rB \PL(n)$ through $\rB \Diff(n) \simeq \rB\rO(n)$. We now quote a nontrivial fact of differential topology, called the \define{Main Theorem of Smoothing Theory}:
\begin{theorem}[{\cite[Essay IV]{MR0645390}}]
\label{smoothingtheory}
  Let $M$ be a piecewise-linear manifold. 
  Then the space of lifts of $\rT_M \colon M \to \rB \PL(n)$ through $\rB \rO(n)$ is homotopy equivalent to the space of smoothings of $M$.
\end{theorem}
An immediate corollary is that a framed PL manifold has a unique  (up to a contractible space) smoothing that is compatible with the framing. In particular, the framed bordism categories built from smooth or from PL manifolds are equivalent. But $\PL(n)$ acts by rotating the framings of $n$-dimensional PL-manifolds, and so acts on the PL version of $\Bord^\fr_n$, and so acts on the space $\iota_0 \cC$. This is the largest group that acts universally on the objects of rigid $(\infty,n)$-categories:
\begin{statement}[{\cite[Remark 2.4.30]{MR2555928}}]\label{thm.PL}
  Assuming the Cobordism Hypothesis, the rotate-the-framing map
  $\PL(n) \to \Aut_{\mathrm{sym\otimes}}(\Bord_n^\fr)$
  is an equivalence when $n\neq 4$. The case $n=4$ is equivalent to the 4-dimensional piecewise-linear Schoenflies conjecture (which remains open).
\end{statement}
The proof of Statement~\ref{thm.PL} has circulated among experts, but is not available in print.

The story in the absence of symmetric monoidal structures is less well-studied. The \define{Tangle Hypothesis} and \define{Cobordism Hypothesis with Singularities} amount to the existence of a graphical calculus for $(\infty,n)$-categories with adjoints, generalizing the graphical calculi described for example in \cite{MR2767048}. The precise details of this graphical calculus have not been worked out in the literature. Roughly, the main ingredients are the following:
\begin{enumerate}
  \item The diagrams in the graphical calculus are drawn on networks of submanifolds of the standard $\bR^n$.
  \item Strata of codimension $k$ are labeled by $k$-morphisms.
  \item Strata are normally framed: if $X \subset \bR^n$ is a codimension-$k$ submanifold, then it comes with a trivialization $\nu \colon \rN_X \cong \bR^k$ of its normal bundle. More generally, substrata of strata are relatively normally framed. This normal framing encodes source and target information (and so must be consistent with the labelings).
  \item Strata are tangentially framed: if $X \subset \bR^n$ is a codimension-$k$ submanifold, then it comes with a trivialization $\tau\colon \rT_X \cong \bR^{n-k}$ of its tangent bundle. This tangential framing encodes the direction of composition internal to the morphism.
  \item The various framing data are compatible. For example, if $X \subset \bR^n$ is a non-sub stratum, it must come equipped with a nullhomotopy of the composite isomorphism
  $$ \bR^n = \rT_{\bR^n} = \rN_X \oplus \rT_X \xrightarrow{\nu\oplus\tau} \bR^k \oplus \bR^{n-k} = \bR^n,$$
  where the left-hand equality simply uses that $X$ is a submanifold of the standard $\bR^n$. For substrata $Y \subset \dots \subset X \subset \bR^n$, there are similar but more complicated compatibility conditions.
\end{enumerate}
When $n=2$, one finds the well known calculus of ``string diagrams.''
The  framing-compatibility is probably the least familiar component of this calculus. Along a 1-dimensional stratum, it is a nullhomotopy of an element of $\rO(2)$. That element is nullhomotopic only when it lives in $\mathrm{SO}(2)$. This forces the normal and tangential framings to determine each other. But there is more to a nullhomotopy than just its existence: an element of $\mathrm{SO}(2)$ has a $\bZ$-torsor of nullhomotopies. Explicitly, the framing compatibility consists of a ``winding number'' carried by each wire. This axiom prevents closed circles; it is what allows the string diagram calculus to apply to non-pivotal 2-categories.

To say that $(\infty,n)$-categories with adjoints admit some graphical calculus is to say that $(\infty,n)$-categories with adjoints are precisely ``interpreters'' for such a graphical calculus. Suppose we are given such an interpreter. Here is another interpreter: precompose your interpreter with some element of $\rO(n)$ acting on all input diagrams. Thus, assuming that there is indeed such a graphical calculus, one finds an action of $\rO(n)$ on $\AdjCat_{(\infty,n)}$. We could act by any element of $\Diff(n)$; since $\rO(n) \hookrightarrow \Diff(n)$ is a homotopy equivalence, this supplies ``the same'' action. But we were not precise about what regularity of manifolds are allowed. 
Statement~\ref{smoothingtheory} and its corollary about unique smoothings for framed PL manifolds says that one might as well work with a PL diagrams, and these are in any case more natural to draw. As such, we find an  action of $\PL(n)$ on $\AdjCat_{(\infty,n)}$. Given Statement~\ref{thm.PL}, we expect:
\begin{conjecture}\label{conj.PL}
  \begin{enumerate}
      \item There is a natural map \label{conj:PL-A}
      \begin{equation*}\label{eqn:PL-action}
          \PL(n) \to \Aut(\AdjCat_{(\infty,n)}).
      \end{equation*}
      \item This map is an equivalence. \label{conj:PL-B}
  \end{enumerate} 
\end{conjecture}
For the remainder of this article, we will assume part (\ref{conj:PL-A}) of Conjecture~\ref{conj.PL}, and we will be motivated by our belief in part (\ref{conj:PL-B}). We emphasize that, not only is the current literature far from a proof of part (\ref{conj:PL-B}), it does not even supply a rigorous construction of this map in part (\ref{conj:PL-A}) except when $n$ is very low:

\begin{example}
 When $n=1$, having adjoints is vacuous, and $\PL(1) = \bZ/2\bZ$ acts by $\cC \mapsto \cC^\op$ as already described. 
\end{example}

\begin{example}\label{eg:O(2)-antiinvolution}
 When $n=2$, $\PL(2) = \rO(2) = \mathrm{SO}(2) \rtimes \bZ/2\bZ = \rB\bZ \rtimes \bZ/2\bZ$. The $\bZ/2\bZ$ subgroup acts by taking $\cC \in \AdjCat_{(\infty,2)}$ to $\cC^{2\op}$, with the opposite composition of $2$-morphisms. Note that there is an equivalence $\cC^{2\op} \simeq \cC^{1\op}$ which is the identity on objects and which acts on 1-morphisms by $f \mapsto f^R$; as such, the assignment $\cC \mapsto \cC^{1\op}$ would produce an equivalent $\bZ/2\bZ$-action on $\AdjCat_{(\infty,n)}$. 
The generator of $\rB \bZ$ acts on $\AdjCat_{(\infty,2)}$ by a natural automorphism of $\id_{\AdjCat_{(\infty,n)}}$. Such a natural automorphism has components: its component at $\cC \in \AdjCat_{(\infty,2)}$ is the autofunctor $\cC \to \cC$ that is the identity on objects and takes a 1-morphism $f$ to its double-right-dual~$f^{RR}$. To fully present an action of $\rB\bZ \rtimes \bZ/2\bZ$ is to present a map out of the space $\rB(\rB\bZ \rtimes \bZ/2\bZ)$, which is the total space of a nontrivial $\bC P^\infty$-bundle over $\bR P^\infty$. 
The fact that $O(2) \ncong O(1) \times SO(2)$ is expressed by the natural isomorphism $((.)^{RR})^{2\op} \cong ((.)^{2\op})^{RR -1}$, using opposites take left to right adjoints.
The higher cells in (a cell model for) this space map to compatibility data between the actions of lower cells.
\end{example}

In analogy with Definition~\ref{defn:volutive}, we declare:
\begin{definition}\label{defn:PL-anti-involutive}
  A $(\infty,n)$-category with adjoints $\cC \in \AdjCat_{(\infty,n)}$ is \define{$\PL(n)$-\wipedge } if it is (equipped with the structure of) a fixed point for the $\PL(n)$-action on $\AdjCat_{(\infty,n)}$. 
\end{definition}

An $(\infty,n-1)$-category, thought of as an $(\infty,n)$-category, never has all adjoints unless it is a groupoid: there is no organic inclusion $\AdjCat_{(\infty,n-1)} \to \AdjCat_{(\infty,n)}$. The other direction, however, does work: if $\cC$ is a  $(\infty,n)$-category with adjoints, then its underlying $(\infty,n-1)$-category $\iota_{n-1}\cC$ has adjoints. We will say that a flagged $(\infty,n)$-category $\cC_0 \to \cC_1 \to \dots \to \cC_n$ \define{has adjoints} if every $\cC_k$ does.

The functor $$\iota_{n-1} \colon \AdjCat_{(\infty,n)} \to \AdjCat_{(\infty,n-1)}$$ cannot be $\PL(n)$-equivariant simply because $\PL(n)$ does not act on the codomain. But it is equivariant for $\PL(n-1)\subset \PL(n)$. Even better: $\iota_{n-1}$
is $(\PL(n-1) \times \PL(1))$-equivariant, where on the domain the action is via the inclusion into $\PL(n)$, and on the codomain the $\PL(1)$-factor acts trivially. The upshot: if $\cC$ is $\PL(n)$-\wipedge, then $\iota_{n-1}\cC$ is $\PL(n-1)$-\wipedge\ and also carries a $\bZ/2\bZ$-action. More generally, the functor $\iota_k\colon \AdjCat_{(\infty,n)}\to \AdjCat_{(\infty,k)}$ is $\PL(k)\times \PL(n-k)$ equivariant, where $\PL(n-k)$ acts on the codomain trivially.  We are led naturally to the following version of Definition~\ref{defn.flagdagncat} for categories with adjoints:

  \setcounter{footnote}{0} 
\begin{definition}\label{defn:PL-dagger}
  A \define{flagged $\PL$-dagger $(\infty,n)$-category}, also called a \define{flagged (fully-)dagger $(\infty,n)$-category with unitary duality}, is a flagged $(\infty,n)$-category with adjoints
  $$ \cC_0 \to \cC_1 \to \dots \to \cC_n$$
  such that each $(\infty,k)$-category $\cC_k$ is $\PL(k)$-\wipedge, and each functor $\cC_k \to \cC_{k+j}$ is $(\PL(k)\times \PL(j))$-\wipedge, with trivial\footnote{More precisely, the $\PL(k) \times \PL(j)$-\wipedgenoun\ on $\cC_k$ is induced by pulling back its $\PL(k)$-\wipedgenoun\ along the projection $\PL(k) \times \PL(j) \to \PL(k)$, and on $\cC_{k+j}$ it is induced by pulling back its $\PL(k+j)$-\wipedgenoun\ along the inclusion $\PL(k) \times \PL(j) \to \PL(k+j)$. } $\PL(j)$-\wipedgenoun\ on $\cC_k$.
\end{definition}

Given a flagged PL-dagger $(\infty,n)$-category 
$$ \cC_0 \to \cC_1 \to \dots \to \cC_n$$
we obtain an underlying flagged fully-dagger $(\infty,n)$-category.
We could ask this fully-dagger to be univalent in the sense of Definition~\ref{defn.dagncat}, which implements the idea that an $(i+1)$-morphism in $\mathcal{C}_i$ is exactly a unitary morphism $\mathcal{C}_{i+1}$.
However, we claim this is not the right univalence condition in general because we want an $(i+2)$-morphism in $\mathcal{C}_i$ to be a $\PL(2)$-unitary morphism in $\mathcal{C}_{i+2}$, and $\PL(2)$ is not isomorphic to $\PL(1) \times \PL(1)$.
Instead, consider the diagram consisting of the  $(\infty,n-l)$-categories $(\iota_{n-l} \cC_{n-j})^{\PL(k_1) \times \dots \PL(k_m)}$, where $(k_1, \dots, k_m)$ is a partition of $j$ for $j \leq l$.
There are two types of maps in the diagram: those forgetting fixed point data corresponding to $(k_1, \dots, k_i+k_{i+1} ,\dots , k_m) \rightsquigarrow (k_1, \dots, k_{i}, k_{i+1}, \dots, k_m)$ and transitions between the flaggings corresponding to $(k_1, \dots, k_{m-1}) \rightsquigarrow (k_1, \dots, k_{m-1}, k_m)$.
Note that if we were to allow the empty partition corresponding to $\iota_{n-l} \mathcal{C}_{n-l} = \mathcal{C}_{n-l}$ we would get a cube. Instead, set $P_{n-l}(\mathcal{C})$ to be the pullback of the diagram built from the nonempty partitions, so that there is a canonical map $\mathcal{C}_{n-l} \to P_{n-l}(\mathcal{C})$. 
For example, when $n=3$ the diagram for $l=3$ is:
\begin{center}
\begin{tikzcd}[sep=small]
 & \mathcal{C}_0\\
	&& P_0(\mathcal{C}) && {\iota_0\mathcal{C}_2^{\PL(2)}} \\
	{\iota_0\mathcal{C}_3^{\PL(3)}} &&& {\iota_0\mathcal{C}_3^{\PL(2)\times\PL(1)}} \\
	&& {\iota_0\mathcal{C}_1^{\PL(1)}} && {\iota_0\mathcal{C}_2^{\PL(1)\times\PL(1)}} \\
	{\iota_0\mathcal{C}_3^{\PL(1)\times\PL(2)}} &&& {\iota_0\mathcal{C}_3^{\PL(1)\times\PL(1)\times\PL(1)}}
	\arrow[from=3-1, to=5-1]
	\arrow[from=5-1, to=5-4]
	\arrow[from=4-3, to=5-1]
	\arrow[from=4-3, to=4-5]
	\arrow[from=4-5, to=5-4]
	\arrow[from=2-5, to=4-5]
	\arrow[from=2-5, to=3-4]
	\arrow[from=2-3, to=3-1]
	\arrow[from=2-3, to=2-5]
	\arrow[from=2-3, to=4-3]
	\arrow[from=3-1, to=3-4,crossing over]
	\arrow[from=3-4, to=5-4, crossing over]
        \arrow[from=1-2, to=2-3, dashed]
\end{tikzcd}
\end{center}
By definition, a morphism in $P_{n-l}(\mathcal{C})$ is a morphism which is $\PL(k)$-unitary for every $k > n-l$. This is precisely what we want for the $(>n-l)$-morphisms in $\mathcal{C}_{n-l}$. Thus we arrive at the following univalence axiom: 

\begin{definition}
  A flagged $\PL$-dagger $(\infty,n)$-category is \define{univalent} if for every $l$, the map $\mathcal{C}_{n-l} \to P_{n-l}(\mathcal{C})$ into the pullback is fully faithful on $> (n-l)$-morphisms. 
  A \define{$\PL$-dagger $(\infty,n)$-category}, also called a \define{dagger $(\infty,n)$-category with unitary duality}, is a univalent flagged one.
\end{definition}

\begin{remark}\label{rem: G-dagger-adjunct}
For good families of groups $G(n)$ related to $\PL(n)$, we expect that there is a definition of $G$-dagger categories with adjoints generalizing Definition~\ref{defn:G-dagger-nonadjunct}. Examples should include, in particular, $\rO(n)$ and $\SO(n)$. We will not try to work out the precise conditions on $G(n)$ or all details of the definition here. 
\end{remark}

Let us turn now to justifying the name ``with unitary duality.'' We will do so by unpacking the notion in the case of bicategories in Examples~\ref{Ex: pivotal bicat} and~\ref{eg:bicat with unitary duals}. Before that, note that the forgetful functor $ \AdjCat_{(\infty,n)}\to \Cat_{(\infty,n)}$ is $(\bZ/2\bZ)^n$ equivariant, where $(\bZ/2\bZ)^n$ acts through the map $(\bZ/2\bZ)^n= \PL(1)^n\to \PL(n)$. 
 Thus to every dagger $(\infty,n)$-category with unitary duality we can assign an underlying (fully-)dagger $(\infty,n)$-category in which the unitary duality is forgotten. So the question is to understand what this extra ``unitary duality'' data looks like.

\begin{example} \label{Ex: pivotal bicat}
We specialize from $(\infty,n)$-categories to bicategories and  look at the $\SO(2)$ subgroup of $\PL(2) = \rO(2) $. What is the data of an $\SO(2)$-\wipedgenoun\ on a bicategory $\cB$? Fixed-point data for an $\SO(2)$($ = \rB \bZ$)-action consists of data assigned to each cell in $\rB \SO(2) = \bC P^\infty$; there is one cell in each even dimension. The $0$-cell selects the bicategory $\cB$. The $2$-cell selects a trivialization $\theta$ of $(-)^{RR}$:
such a trivialization unpacks to a family of 1-isomorphisms $\theta_b\colon b\to b $ for all $b\in \mathcal{B}$ and 2-isomorphisms $ f^{RR} \circ \theta_{b_1} \to \theta_{b_2} \circ f$ for $1$-morphisms $f\colon b_1 \to b_2$ satisfying various coherence conditions. The $4$-cell selects a quadratic equation that $\theta$ must solve; the reader is encouraged to work out this equation as an exercise. In the special case of bicategories, this is all the necessary data: the higher cells in $\bC P^\infty$ admit unique assignments, because the space of bicategories with adjoints is a homotopy 3-type. 

The group $\SO(1)$ is trivial, and so to enhance an $\SO(2)$-\wipedge\ bicategory $\cB$ to an $\SO(2)$-dagger bicategory, one needs only to supply the data of an $\SO(2)$-\wipedge~ essentially surjective functor $\cB_0\to \cB$ such that $\cB_0\to (\iota_0 \cB)^{\SO(2)}$ is 1-fully faithful.
The objects $(\iota_0 \cB)^{\SO(2)}$ are given by pairs consisting of of an object $b\in \cB$ and a 2-isomorphism $\omega_b\colon \theta_b \to \id_b$. Similar to Section~\ref{sec.defbicat}, form a new bicategory $\cB'$ whose objects are those of $\cB_0$ and morphisms are those in $\cB$. This bicategory $\cB'$ comes with a trivialization of the double dual functor $(-)^{RR}$ which is the identity on objects (and which solves a quadratic equation). In other words, $\cB'$ is a \define{pivotal bicategory} as defined for example in~\cite[Definition 2.1]{pivotalstringnets}, also called an \define{even-handed bicategory}~\cite{bartlett}. Like traditionally-defined dagger $1$-categories, pivotal bicategories are ``evil'' in the sense that pivotal structures do not transport across bicategorical equivalences. 
Remark~\ref{rem: G-dagger-adjunct} suggests
a coherent version of ``pivotal bicategory'': they are the $\SO(2)$-dagger bicategories.
\end{example}

To finish the discussion of the name ``unitary duality,'' we now restore the reflection. 
\begin{example}\label{eg:bicat with unitary duals} 
  As mentioned in Example~\ref{eg:O(2)-antiinvolution}, to discuss actions by, and fixed points for, $\rO(2) = \rB \bZ \rtimes \bZ/2\bZ$, one should assign data to the cells in a cell model for $\rB \rO(2)$, which is a $\bC P^\infty$-bundle over $\bR P^\infty$. In particular, $\rB \rO(2)$ has a cell model with one cell of dimension $2s+r$ for each pair $s,r \in \bN$. The cells indexed $(s,0)$ supply the restriction of the data along $\SO(2) \to \rO(2)$, and the cells indexed $(0,r)$ supply the restriction of the data along $(\bZ/2\bZ)_2 \subset \rO(2)$.
  We see that an $\rO(2)$-volution consists of a $(\bZ/2\bZ)_2$-volution $\psi_2\colon \cB \to \cB^{2\op}$, a trivialization $f^{RR} \to \theta_{b_2} \circ f \circ \theta_{b_1}^{-1}$ as in Example \ref{Ex: pivotal bicat}, together with a natural modification equating the two ways to identify $\psi_2(f^{RR}) \cong \psi_2(f)^{LL}$ with $\psi_2(f)$,
   satisfying some conditions.
  This in particular includes a $(\bZ/2\bZ)_2$-fixed point datum $\psi_2(\theta_b)) \cong \theta_{\psi_2(b_1)}^{-1}$.
  We see that an $\rO(2)$-dagger structure consists of
  \begin{enumerate}
    \item an $\rO(2)$-equivariant essentially surjective functor $\cB_0 \to \cB$ from a $2$-groupoid such that $\cB_0 \to (\iota_0 \cB)^{\rO(2)}$ is fully faithful. 
    The $\rO(2)$-fixed points combine the trivializations $\omega$ of Example \ref{Ex: pivotal bicat} from $\SO(2)$ with the $\rO(1)$-fixed point data $h_b^2$ such that $\omega$ is compatible with $\psi_2(\theta_b)) \cong \theta_{\psi_2(b_1)}^{-1}$ and the $2$-isomorphism $f^{RR} \to \theta_{b_2} \circ f \circ \theta_{b_1}^{-1}$ preserves $(\bZ/2\bZ)_2$-fixed point data. 
      \item As in Section \ref{sec.defbicat}, the $(2,1)$-category $\cB_1$ and the functors $\cB_1 \to \cB$ and $\cB_0 \to \cB_1$ are fixed after specifying $(\bZ/2\bZ)_2$-fixed point data on $1$-morphisms of $\cB$. 
    \item The functors $\cB_0 \to \cB_1$ and $\cB_1 \to \cB$ are still required to be $(\bZ/2\bZ)^2$-equivariant. As explained in Section \ref{sec.defbicat}, the $(\bZ/2\bZ)_2$ simply specifies the agreement of fixed point data on $2$-morphisms. 
    For the $(\bZ/2\bZ)_1$-equivariance, first note that the underlying $(\bZ/2\bZ)_1$-volution of the $\rO(2)$-volution can be described by $\psi_1 = (.)^R \circ \psi_2$, where we made an arbitrary choice of the right adjoint to trivialize the $180^\circ$ rotation in $\SO(2)$.
    Using the pivotal structure and the fact that $\psi_2(f^R) \cong \psi_2(f)^L$, we obtain data specifying that $\psi_2$ and $(.)^R$ commute.
    Note that the $(\bZ/2\bZ)_1$-fixed point data induced by the $(\bZ/2\bZ)_2$-fixed point data $h^2_b$ is $(h^2_b)^{R}$. 
    The $(\bZ/2\bZ)_1$-equivariance data of $\cB_1 \to (\iota_1 \cB)^{(\bZ/2\bZ)_2}$ will ensure that if $h_f$ is a $(\bZ/2\bZ)_2$-fixed point data on a $1$-morphism, then there is a canonical fixed point data $h_{f^R}$ for $f^R$, which makes the right adjoint into a top dagger functor.
  \end{enumerate}

  When the bicategory $\cB$ is the one-object delooping of a rigid monoidal 1-category $\cC$, and $\cB_0$ is declared to consist exactly of this object, this should recover the notion of a dagger category with a unitary dual functor~\cite[Definition on page 53 part (2)]{MR2767048}, see also~\cite{MR4133163}.
\end{example}

As in Definition~\ref{defn:coflagged}, one can also define \define{coflagged} $\PL$-dagger $(\infty,n)$-categories by keeping the univalence axiom Definition~\ref{defn:PL-dagger} but dropping the essential surjectivity (i.e.\ not requiring that $\cC_0 \to \dots \to \cC_n$ be a flagging). 
Statement~\ref{statement:completioncoflagging} generalizes: every coflagged $\PL$-dagger category can be completed to a $\PL$-dagger category by replacing each $\cC_k$ with the full image of $\cC_{k-1}$; this is the right adjoint to the forgetful functor from $\PL$-dagger to coflagged $\PL$-dagger categories. The further forgetful functor from coflagged $\PL$-dagger to $\PL$-anti-involutive categories that remembers only the top level also has a right adjoint which assigns to a $\PL$-anti-involutive $(\infty,n)$-category $\cC$ the diagram
\begin{equation}\label{Eq: herm completion}
(\iota_0 \cC)^{\PL(n)}\to (\iota_1 \cC)^{\PL(n-1)} \to \dots \to  \cC   .
\end{equation}

\section{Bordism categories and reflection-positive topological quantum field theories}\label{sec:bord}
For some of us, our interest in higher dagger categories stems from the role that we expect them to play in the study of (topological) quantum field theories.
 For any space $X\to \rB\rO(n)$ over $\rB \rO(n)$ there is a rigid symmetric monoidal $(\infty,n)$-category $\Bord_n^{X}$ of bordisms $M$ with a lift of their tangent bundle $\rT_M \colon M \to \rB\rO(n)$ through $X$. (The construction of $\Bord_n^X$ is outlined in \cite{MR2555928} and fully implemented by one of us in \cite{MR3924174}.) For example, a framing is a trivialization of $\rT_M$, i.e.\ a lift through $\ast \to \rB\rO(n)$, and $\Bord_n^\ast = \Bord_n^\fr$.
  Let $\cC$ be a symmetric monoidal $(\infty,n)$-category. An \define{$X$-structured $n$-dimensional fully extended topological quantum field theory with values in $\cC$} is a symmetric monoidal functor $\mathcal{Z}\colon \Bord_n^{X}\to \cC$. Since the bordism category is rigid, every topological quantum field theory will factor through the maximal rigid sub-$(\infty,n)$-category of $\cC$. Hence we will from now on assume without loss of generality that $\cC$ is rigid. 
Recall that the Cobordism Hypothesis allows us to identify fully extended framed topological field theories with $\iota_0\cC$, which hence carries a $\PL(n)$-action. 
This was the starting point for the discussion in the previous section. The Cobordism Hypothesis also provides a description of the space of $X$-structured fully extended topological quantum field theories as the space of $X$-fixed points $(\iota_0 \cC)^X$. 

As an example, let us focus on unoriented smooth $n$-dimensional topological quantum field theories valued in $\cC$. The space of such theories is $(\iota_0\cC)^{\rO(n)}$. The $1$-morphisms in $\cC$ are the objects of a rigid symmetric monoidal $(\infty,n-1)$-category $\iota_{n-1}\cC^\to$ \cite[Section 7]{MR3590516}, and thus carry an action by $\rO(n-1)$; the maps that assign to a 1-morphism its source and target and that  compose 1-morphisms  are $\rO(n-1)$-equivariant, and compile into an action of $\rO(n-1)$ on the whole $(\infty,1)$-category $\iota_1 \cC$. This action is meaningful from the perspective of topological quantum field theory: there is an $\rO(n-1)$-action on the space of 1-morphisms between the objects underlying two different unoriented topological field theories, and the Stratified Cobordism Hypothesis says that fixed points for this $\rO(n-1)$-action classify unoriented codimension-$1$ defects between the two field theories. Similarly, $2$-morphisms carry an induced $\rO(n-2)$-action whose fixed points classified unoriented codimension-$2$ defects, and so on.
 In summary, unoriented smooth defects assemble into the following diagram:
\begin{equation} \label{eqn:unorienteddefects}
(\iota_0 \cC)^{\rO(n)}\to (\iota_1 \cC)^{\rO(n-1)} \to \dots \to  \cC  .
\end{equation}
There is nothing special in~\eqref{eqn:unorienteddefects}  about unoriented smooth theories and the groups $\rO(k)$: any good sequence of groups $G(k)$ would work (compare Remark~\ref{rem: G-dagger-adjunct}). For example, replacing the $\rO(k)$'s in~\eqref{eqn:unorienteddefects} with $\SO(k)$'s would organize the defects between oriented smooth field theories as analyzed for example in~\cite{MR2742426,MR4065261}; replacing the $\rO(k)$'s with $\PL(k)$'s would organize the defects between unoriented piecewise linear field theories.

The diagram~\eqref{eqn:unorienteddefects} is suspiciously close to~\eqref{Eq: herm completion} but they are a priori different: the action of $\rO(k)$ on $\iota_{n-k}\cC$ in~\eqref{eqn:unorienteddefects} comes from the Cobordism Hypothesis and uses the symmetric monoidal structure on $\cC$, whereas in~\eqref{Eq: herm completion} we envisioned selecting a $\PL(n)$-\wipedgenoun\ on $\cC$, i.e.\ fixed-point data for the action of $\PL(n)$ on the graphical calculus for $\cC$. We expect the relation to be the following. 
Any symmetric monoidal $(\infty,n)$-category $\cC$ determines a symmetric monoidal $(\infty,n+1)$-category $\rB \cC$ with one object $\bullet$ and $\operatorname{End}_{\rB \cC}(\bullet) = \cC$; iterating this supplies $(\infty,n+m)$-categories $\rB^m\cC$ for every $m \in \bN$. 
If $\cC$ is rigid, then $\rB^m\cC$ will have all adjoints. The graphical calculi for $\cC$ and for $\rB^m \cC$ are compatible via the embedding of $\bR^n$ into $\bR^{n+m}$ as the last $n$ coordinates. Hence we can think of the graphical calculus for a rigid symmetric monoidal $(\infty,n)$-category as taking place in $\bR^{\infty}$. The space of embeddings (of any finite-dimensional object) into $\bR^\infty$ is contractible. This contractibility selects a canonical trivialization of the $\PL(n)$-action predicted in Conjecture~\ref{conj.PL} for any rigid symmetric monoidal $(\infty,n)$-category. In other words, if we let $\SymRigidCat_{(\infty,n)}$ denote the $(\infty,1)$-category of rigid symmetric monoidal $(\infty,n)$-categories: 

\begin{statement}
    The forgetful functor $$\SymRigidCat_{(\infty,n)} \to \AdjCat_{(\infty,n)}$$ factors through the $\PL(n)$-fixed points $\AdjCat_{(\infty,n)}^{\PL(n)} \to \AdjCat_{(\infty,n)}$. 

In particular, a rigid symmetric monoidal structure on a $(\infty,n)$-category selects a canonical $\PL(n)$-\wipedge\ structure on its underlying $(\infty,n)$-category with adjoints. Moreover, any symmetric monoidal functor between rigid symmetric monoidal $(\infty,n)$-categories will automatically intertwine these canonical $\PL(n)$-\wipedge\ structures.
\end{statement}

Given a rigid symmetric monoidal $(\infty,n)$-category $\cC$, we expect that the induced $\PL(n)$-action on $\iota_0\cC$ agrees with the one coming from the Cobordism Hypothesis. 

Other symmetric monoidal $\PL(n)$-\wipedge\ structures  on $\cC$ are given by twisting the canonical one by a symmetric monoidal $\PL(n)$-action on $\cC$: symmetric monoidal $\PL(n)$-\wipedge\ structures  form a trivialized torsor over $\PL(n)$-actions. 

There is a straightforward way to define symmetric monoidal $\PL(n)$-dagger categories by simply requiring the \wipedgenoun, the flagging by $\cC_i$'s, and the trivialization data, to all be symmetric monoidal.
Interestingly, if $\cC$ is in addition rigid
the sequence of deloopings $\cC, \rB\cC, \rB^2\cC,\dots$ allows more. Suppose that 
\begin{equation} \label{Aflagging} \cC_0 \to \cC_1 \to \dots \to \cC_n = \cC\end{equation}
is a  symmetric monoidal flagged $(\infty,n)$-category. Then the deloopings assemble into a flagged $(\infty,n+1)$-category by selecting the unit object at the bottom:
\begin{equation} \label{Bflagging} \{\bullet\} \to \rB\cC_0 \to \rB \cC_1 \to \dots \rB\cC_n = \rB \cC.\end{equation}
If $\cC_k$ is rigid as a symmetric monoidal $(\infty,k)$-category, then $\rB \cC_k$ is also rigid, now as a symmetric monoidal $(\infty,k+1)$-category. 

Given such a structure, just like it is natural to ask for a symmetric monoidal flagged $\PL(n)$-dagger structure on a chosen flagging~\eqref{Aflagging}, it is natural to ask to give $\rB\cC$ a flagged $\PL(n+1)$-dagger structure with chosen flagging ~\eqref{Bflagging}. Iterating, it is natural to ask to give $\rB^m\cC$ a compatible symmetric monoidal $\PL(n+m)$-dagger structure. We call the structure just sketched a \define{\RENAMEME\ with underlying $(\infty,n)$-category~$\cC$}.

After untwisting by the canonical \wipedge\ structures, a symmetric monoidal $\PL(n+m)$-dagger structure on $\rB^m \cC$ unpacks to a symmetric monoidal action of $\PL(k+m)$ on $\cC_k$ for all $k \leq n$, and compatible $(\PL(k+m) \times \PL(j))$-equivariant maps $\cC_k \to \cC_{k+j}$ for all $j + k \leq n$. Univalence is just about the $\PL(j)$-fixed points: it is independent of $m$.  A more careful treatment and definition is beyond the scope of this short article. 

\begin{example}
To illustrate the difference between the two definitions above, we look at dagger 1-categories. 
A (strict) symmetric monoidal dagger 1-category is a symmetric monoidal category equipped with a symmetric monoidal anti-involution which is the identity on objects.
From Example~\ref{eg:bicat with unitary duals} we learn that a \RENAMEME\ with underlying 1-category is in addition equipped with a unitary dual functor. There are no additional structures corresponding to higher deloopings, since they would correspond to higher morphisms.   
\end{example}

We will now construct a dagger structure on certain bordism categories $\Bord_n^X$.  
Let us recall some more details about these $(\infty,n)$-categories. They depend on a choice of tangential structure $X \to \rB \PL(n)$---for example, by smoothing theory (Theorem~\ref{smoothingtheory}), the smooth unoriented bordism category $\Bord_n^{\rO(n)}$ corresponds to the tangential structure $\rB \rO(n) \to \rB \PL(n)$. Given $k \leq n$, set $X(k) \to \PL(k)$ to be the pullback of $X \to \PL(n)$ along the standard inclusion $\rB\PL(k) \to \rB\PL(n)$. A $k$-morphism in $\Bord_n^X$ is by definition a $k$-dimensional cobordism $M^k$ equipped with a lift of its tangent bundle $\rT_M\colon M \to \rB\PL(k)$ through $X(k)$; we will refer to such a lift simply as an \define{$X$-structure}.

Just as the tangential structure $X(k)$ on $k$-morphisms is pulled back from an $n$-dimensional tangential structure, it often happens that the input $n$-dimensional tangential structure is itself pulled back from a higher-dimensional tangential structure:
\begin{definition}\label{defn:stabletangentialstructure}
    A \define{stabilization} of an $n$-dimensional tangential structure $X(n) \to \rB\PL(n)$ is a map $X \to \rB\PL = \injlim_{n\to\infty} \rB\PL(n)$ and an equivalence $X(n) \cong X \times_{\rB\PL} \rB \PL(n)$. A tangential structure is \define{stable} when it is equipped with a stabilization.
\end{definition}
We emphasize that stability is structure, not just a property.

\begin{example}
Define a \define{stable smooth structure} on a PL-manifold $M$ to be a lift of its stabilized tangent bundle through $\rB\rO \to \rB\PL$. In other words, for an $n$-manifold, stable smoothness corresponds to the tangential structure $\rB\rO \times_{\rB\PL} \rB\PL(n)$. Stable smoothness is by construction a stable tangential structure.

Stable smoothness is weaker than smoothness: the map $\rB\rO(n) \to \rB\rO \times_{\rB\PL} \rB\PL(n)$ is not a homotopy equivalence. Indeed, actual-smoothness is not a stable tangential structure. Write $\Bord_n^{\rO(n)}$ for the (unoriented) actually-smooth bordism $(\infty,n)$-category and $\Bord_n^\rO$ for the stably-smooth version. The map $\Bord_n^{\rO(n)} \to \Bord_n^\rO$ is not an equivalence of $(\infty,n)$-categories.

That said, a careful analysis of smoothing theory shows that the map $\rB\rO(n) \to \rB\rO \times_{\rB\PL} \rB\PL(n)$ is $n$-connected; in particular, this follows from \cite[Lecture 21, Theorem 1]{LurieSmoothing}. 
In terms of bordism categories, this implies that $\Bord_n^{\rO(n)} \to \Bord_n^\rO$ becomes an equivalence after quotienting to weak $n$-, a.k.a. $(n,n)$-, categories. The upshot is that for TQFTs valued in a weak $n$-category, the notions of smooth and stably-smooth TQFT agree.
\end{example}

The description of $k$-morphisms in $\Bord_n^X$ as $X$-structured $k$-dimensional cobordisms does not precisely present the actual space of $k$-morphisms in $\Bord_n^X$. Indeed, the natural space of $k$-dimensional $X$-structured bordisms has as its equivalences the $X$-structured diffeomorphisms. But the categorical equivalences in $\Bord_n^X$ between $(<n)$-dimensional morphisms also include $X$-structured h-cobordisms, and in high dimensions not every h-cobordism comes from a diffeomorphism. Rather, this description presents a flagged $(\infty,n)$-category of bordisms as
\begin{equation} \label{eqn:bordflagged}
     \Bord_0^X \to \Bord_1^X \to \dots \to \Bord_n^X.
\end{equation}
The main result of this section is:

\begin{statement}\label{statement:bordismisdagger}
  Suppose that $X \to \rB\PL$ is a stable tangential structure. Then the flagged $X$-structured bordism category~\eqref{eqn:bordflagged} is naturally a flagged $\PL(n)$-dagger $(\infty,n)$-category.
\end{statement}

\begin{proof}[Construction]
  To build a dagger  on $\Bord_{n}^X$, we must first provide a $\PL(n)$-\wipedgenoun. 
  We will do this by twisting the canonical $\PL(n)$-\wipedgenoun\ by a carefully-selected action of $\PL(n)$ on $\Bord_n^X$. To build such an action, consider the functor $$\cat{Space}_{/\rB\PL} \to \cat{Space}_{/\rB\PL(n)} \to \cat{SymCat}_{(\infty,n)} \hspace{1cm} X \mapsto \Bord_{n}^X$$ which first pulls back a stable tangential structure $\{X\to \rB \PL\} \in \cat{Space}_{/\rB \PL}$ to its unstable variant $X(n):= X\times_{\rB \PL} \rB \PL(n)$ and then assigns the corresponding bordism category. 
  To build an action of $\PL(n)$ on the resulting bordism category, it therefore suffices to specify a functor 
  $$
  \rB \PL(n) \to \cat{Space}_{/\rB \PL}
  ,$$ i.e. an action of $\PL(n)$ on the bundle $X \to \rB\PL$. Equivalently, we want to choose a bundle $X' \to \rB\PL \times \rB\PL(n)$ whose pullback along $\rB\PL = \rB\PL \times \{\pt\} \to \rB\PL \times \rB\PL(n)$ is $X$. 
  
  Recall that $\rB\PL$ is a coherently-associative monoid under $\oplus$, and that in homotopy theory a monoid is a group as soon as its $\pi_0$ is a group. Thus $\rB\PL$ comes with an operation $\text{``}{\ominus}\text{''}\colon \rB\PL \to \rB \PL$, canonical up a contractible space of choices, that inverts with respect to the group operation $\oplus$. With this operation in hand, we choose to set $X'$ to be the pullback of $X \to \rB \PL$ along the composition
  $$ \rB\PL \times \rB\PL(n) \to \rB\PL \times \rB \PL \xrightarrow{\id \times \ominus} \rB \PL \times \rB \PL \xrightarrow{\oplus} \rB \PL.$$

  These actions are compatible: the inclusion $\Bord_k^X \to \Bord_n^X$ intertwines the induced $\PL(k)$-\wipedgenoun s. What remains is to trivialize a certain induced action of $\PL(n-k)$ on $\Bord_{k}^X$. Our $\PL(n)$-\wipedgenoun\ on $\Bord_n^X$ was built out of two pieces: the canonical \wipedgenoun\ together with our choice $X'$. By the same token, the action of $\PL(n-k)$ on $\Bord_{k}^X$ that we care about has two pieces. One piece is the restriction along $\PL(n-k) \to \PL(n)$ of the action on $X \to \rB \PL$. The other piece comes from restricting the canonical \wipedgenoun\ along $\PL(n-k) \to \PL(n)$.
  Combining these pieces and writing the problem in terms of bundles, one finds that what needs trivializing is the restriction of $X$ along the composition
  $$ \rB\PL(n-k) \xrightarrow{\operatorname{diag}} \rB \PL(n-k) \times \rB \PL(n-k) \to \rB\PL \times \rB \PL \xrightarrow{\id \times \ominus} \rB \PL \times \rB \PL \xrightarrow{\oplus} \rB \PL.$$
  But this composition factors through
  $$ \rB\PL \xrightarrow{\operatorname{diag}} \rB\PL \times \rB \PL \xrightarrow{\id \times \ominus} \rB \PL \times \rB \PL \xrightarrow{\oplus} \rB \PL,$$
  which is trivial by the definition of $\ominus$.
\end{proof}

\begin{remark}\label{Rem: smooth instead of PL}
Both the notion of a dagger category with unitary duals and 
Statement~\ref{statement:bordismisdagger} are based on $\PL$ geometry. There is an analogous story replacing $\PL(n)$ with $\rO(n)$ and $\PL$-geometry with smooth geometry. The smooth bordism category with smoothly stable tangential structure has a natural $\rO$-dagger structure.     
\end{remark}

\begin{example}\label{eg:bordism}
Let us restrict to unextended bordism categories. Whereas the very definition of $\PL(n)$-dagger $(\infty,n)$-category, and hence the content of Statement~\ref{statement:bordismisdagger}, required Conjecture~\ref{conj.PL}, the definition in the unextended case and the construction described below are fully rigorous.

Because smooth bordism categories are more familiar than their piecewise-linear counterparts, we will focus on that case.
For each $X \to \rB\rO(n)$ there is a smooth $X$-structured unextended bordism $(\infty,1)$-category $\Bord_{n,n-1}^X$ whose objects are closed $(n-1)$-dimensional $X$-structured manifolds and whose morphisms are $X$-structured $n$-dimensional bordisms. This category is naturally flagged: the na\"ive equivalences between objects in $\Bord_{n,n-1}^X$ are the $X$-structured diffeomorphisms, but the categorical equivalences also include $X$-structured h-cobordisms; writing $\operatorname{Man}_{n-1}^X$ for the space of $X$-structured closed $(n-1)$-manifolds and diffeomorphisms between them, we find a  flagged $(\infty,1)$-category
$$     \operatorname{Man}_{n-1}^X \to \Bord_{n,n-1}^X$$
  which is often not univalent. 

Since $\Bord_{n,n-1}^X$ is rigid symmetric monoidal, it has a canonical anti-involution (which is to say a \wipedgenoun\ for the group $\PL(1) = \bZ/2\bZ$). It is precisely the operation $(-)^\vee$ that takes duals.
This canonical anti-involution typically does not extend to a dagger structure: for example, when $X = \SO$, the dual of an object is its orientation-reversal, and most oriented manifolds simply are not orientation-diffeomorphic to their orientation reversals. This is why we needed to find an appropriate twist by a symmetric monoidal action.
  
  Suppose that our tangential structure is stable, i.e.\ (in the smooth case) pulled back from $\rB\rO$. We will twist the canonical anti-involution on $\Bord_{n,n-1}^X$ by the involution on $X \to \rB\rO$ built by restricting along the composition
  \begin{equation} \label{eqn:smoothunextended}
   \rB\rO \times \rB(\bZ/2\bZ) = \rB\rO \times \rB\rO(1) \to \rB \rO \times \rB \rO \xrightarrow{\id \times \ominus} \rB \rO \times \rB \rO \xrightarrow{\oplus} \rB \rO. \end{equation}
   
  To finish the construction of a dagger structure, we take the corresponding twisted anti-involution, restrict it to the (flagged) groupoid $\operatorname{Man}_{n-1}^X$ of objects where it becomes an involution, and trivialize it. This involution combines duality with the chosen  involution on $X \to \rB\rO$. 
  
  In general, given an $n$-dimensional tangential structure $X \to \rB\rO(n)$, the duality involution on $\operatorname{Man}_{n-1}^X$ is the one built from restricting $X$ along $\rB\rO(n-1) \times \rB(\bZ/2\bZ) = \rB\rO(n-1) \times \rB\rO(1) \to \rB \rO(n)$. Together with the twisting, we find that we win if the total composition
  $$ \rB\rO(n-1) \times \rB(\bZ/2\bZ) \xrightarrow{\id \times \operatorname{diag}} \rB\rO(n-1) \times \rB(\bZ/2\bZ) \times \rB(\bZ/2\bZ) \to \rB\rO(n) \times \rB(\bZ/2\bZ) \xrightarrow{\eqref{eqn:smoothunextended}} \rB \rO$$
  agrees with the canonical map $\rB\rO(n-1) \to \rB\rO$. And it does by virtue of $\ominus$.
\end{example}

  \begin{statement}
The twist defined in (\ref{eqn:smoothunextended}) is the same twist used to form $\hat{H}^{(1)}_n$ in \cite[Appendix E]{MR4268163}, and hence our dagger structure agrees with the dagger structure implicit in \cite{MR4268163}.
\end{statement}

Symmetric monoidal dagger structures on bordism categories are relevant for the definition of reflection positivity~\cite{MR1001453, MR2387729, MR3674995, MR4268163, luukthesis}. 
As a motivating example,
consider the anti-involution on the symmetric monoidal category of finite dimensional super vector spaces $\sVect$ given by twisting the canonical one by the $\bZ/2\bZ$-action corresponding to complex conjugation.  
The category of fixed points $(\iota_0 \sVect)^{\bZ/2}$ describes Hermitian super vector spaces and unitary maps between them. 
We denote by $\sHilb^{u} \subset (\iota_0 \sVect)^{\bZ/2}$ its full subcategory on the super Hilbert spaces and by $\sHilb$ the coherent dagger 1-category
$$
 \sHilb^{u} \to \sVect.
$$
An $X$-structured unextended topological quantum field theory $Z \colon \Bord^X_{n,n-1} \to \sVect$ has a \define{reflection structure} when it is a functor of anti-involutive, or equivalently $\bZ/2\bZ$-equivariant, categories.  A field theory with reflection structure is called \define{reflection positive} if it induces a dagger functor $\Bord^X_{n,n-1}\to \sHilb$: if the Hermitian structures induced by Example~\ref{eg:bordism} on the images of $Z$ are all positive definite. Reflection positivity encodes physical unitarity. Thus we propose the following tentative definition:
\begin{proposal}
\label{defn:unitaryTQFT}
    Let $\cC$ be a rigid symmetric monoidal $\PL(n)$-dagger $(\infty,n)$-category. A \define{unitary $X$-structured extended topological quantum field theory} valued in $\cC$ is a functor of symmetric monoidal $\PL(n)$-dagger categories $\Bord_n^X \to \cC$.
\end{proposal}

\begin{remark}
Proposal~\ref{defn:unitaryTQFT} is to be interpreted loosely as the type of definition we expect to be correct. However, there are many details to be pinned down. For instance, we expect the dagger structure on $\Bord_n^X$ from Statement~\ref{statement:bordismisdagger} to be part of a \RENAMEME\  and a unitary extended topological field theory might be required to respect this structure.  
In addition, our proposal prefers the $\PL$-bordism category, but the smooth version from Remark~\ref{Rem: smooth instead of PL} might be more appropriate in certain applications.  
\end{remark}

Fully defining unitary extended topological quantum field theory furthermore requires settling on good target categories $\cC$. Without unitarity, good target categories should feel like higher categories of ``higher vector spaces'' or ``higher super vector spaces''; with unitarity, one should instead expect $\cC$ to organize the ``higher (super) Hilbert spaces.'' In \cite{2410.05120}, some of us construct a symmetric monoidal $(\infty,3)$-category $\cat{Hilb}_3$ of ``finite-dimensional 3-Hilbert spaces.'' It is rigid by construction, and moreover it is expected to carry an organic $\PL(3)$-dagger structure. Together with Definition~\ref{defn:unitaryTQFT} and Statement~\ref{statement:bordismisdagger}, one arrives at a definition of $3$-dimensional bosonic fully-extended unitary topological quantum field theories for any stable tangential structure.  The construction of $\cat{Hilb}_3$ in~\cite{2410.05120} starts with the usual dagger $1$-category $\cat{Hilb}$ of finite-dimensional Hilbert spaces and repeatedly applies a certain manifestly-dagger delooping procedure. 
This delooping procedure is expected to extend to even higher categories as well. 

Another candidate for $\cat{Hilb}_3$ weakening the finite dimensionality conditions would be an appropriate Morita $3$-category of Bicommutant Categories. Bicommutant categories as first introduced in \cite{MR3747830}; their Morita theory is under development. 
Similar to the 3-category $\cat{TensCat}$, the category of Bicommutant categories, bicommutant category bimodules (W*-categories), equivariant functors and natural transformations, called $\cat{BicomCat}$ could serve as target for unitary quantum field theories, more general than just TQFTs, for instance chiral conformal field theories. It is expected that this category is a strictly fully-dagger, with involutions given by the categorical $\mathrm{}^\mathrm{op}$ operation on objects and 1-morphisms and using adjoint functors and natural transformations on the top two levels.

Further analysis of Definition~\ref{defn:unitaryTQFT}---construction of examples, a Unitary Cobordism Hypothesis, etc.\ ---will be the subject of future work.

\subsection*{Acknowledgments}
This article evolved out of the June 2023 online workshop ``Dagger Higher Categories.'' We gratefully thank the other workshop participants who chose not to join in the writing on this article: Bruce Bartlett, 
Andr\'e Henriques,
Chris Heunen,
Peter Selinger,
and
Dominic Verdon. We also thank David Ayala and Jan Steinebrunner for inspiring discussions on related topics.

The authors acknowledge the following grant support:
\begin{center}  
\footnotesize{
\begin{tabular}{cc} \\[-18pt]
GF, BH, DP & NSF DMS 2154389  
\\
CK & NSF GRFP 2141064 \\
TJF & NSERC RGPIN-2021-02424\\
TJF, LM, CS, LS & Simons Collaboration on Global Categorical Symmetries (Simons Foundation grants 888996 and 1013836)\\
LS & Atlantic Association for Research in the Mathematical Sciences
\\
DR & Deutsche Forschungsgemeinschaft (DFG) – 493608176
\\
CS & Deutsche Forschungsgemeinschaft (DFG) – SFB 1085 Higher Invariants \\
N  & James Buckee Scholarship, Merton College, Oxford\\
\end{tabular}
}
\end{center}
\footnotesize{Research at Perimeter Institute is supported in part by the Government of Canada through the Department of Innovation, Science and Economic Development and by the Province of Ontario through the Ministry of Colleges and Universities.}
\bibliography{biblio}
\bibliographystyle{alpha}

\end{document}